\documentclass[journal]{IEEEtran}

\usepackage{cite}
\usepackage{amssymb,amsmath, amsfonts, amsthm}
\usepackage[caption=false,font=footnotesize]{subfig}
\usepackage{tikz,pgfplots}
\usepackage{pgfplots}
\usetikzlibrary{shapes.geometric, decorations.pathmorphing, positioning, arrows, arrows.meta, patterns, shapes, calc, fit, backgrounds}
\graphicspath{{./figures/TAP_TDCFIE/}}
\DeclareGraphicsExtensions{.pdf,.jpeg,.png,.eps}

\newcommand{\R}{\mathbb R}

\newcommand{\G}{\mathbb G}
\newcommand{\M}{\mathbb M}

\newcommand{\Z}{\mathbb Z}
\newcommand{\X}{\mathbb X}
\newcommand{\T}{\mathbb T}
\newcommand{\Pbb}{\mathbb P}

\newcommand{\MM}{\mathcal M}

\newcommand{\TT}{\mathcal T}
\newcommand{\OO}{\mathcal O}

\newcommand{\KK}{\mathcal K}
\newcommand{\ZZ}{\mathcal Z}
\newcommand{\XX}{\mathcal X}

\newcommand{\yrm}{\mathrm{y}}
\newcommand{\jrm}{\mathrm{j}}

\newcommand{\erm}{\mathrm{e}}
\newcommand{\hrm}{\mathrm{h}}

\DeclareMathOperator{\xxs}{\mathbf x}

\DeclareMathOperator{\yys}{\mathbf y}
\DeclareMathOperator{\rrs}{\mathbf r}
\DeclareMathOperator{\jjs}{\mathbf j}

\DeclareMathOperator{\HHs}{\mathbf H}

\DeclareMathOperator{\XXs}{\mathbf X}

\DeclareMathOperator{\PPo}{\mathbf P}

\DeclareMathOperator{\ZZs}{\mathbf Z}

\newcommand{\paren}[1]{\left({#1}\right)}

\newcommand{\abs}[1]{\left\vert{#1}\right\vert}

\newcommand{\inprod}[1]{\left\langle{#1}\right\rangle}

\newcommand{\ovl}[1]{\overline{#1}}

\newcommand{\transpose}{\mathsf{T}}

\newcommand{\bm}[1]{\boldsymbol{#1}}

\newcommand{\curlt}{\textbf{\textup{curl}}}
\newcommand{\divt}{\textup{div}}
\newcommand{\gradt}{\textbf{grad}}

\newcommand{\pa}{\partial}

\newcommand{\Gm}{\Gamma}

\newcommand{\Sm}{\Sigma}
\newcommand{\vphi}{\varphi}

\newcommand{\rb}{\bm{r}}
\newcommand{\hb}{\bm{h}}

\newcommand{\eb}{\bm{e}}

\newcommand{\lb}{\bm{l}}
\newcommand{\fb}{\bm{f}}
\newcommand{\gb}{\bm{g}}
\newcommand{\kb}{\bm{k}}
\newcommand{\jb}{\bm{j}}
\newcommand{\ib}{\bm{i}}

\newcommand{\pb}{\bm{p}}
\newcommand{\ssb}{\bm{s}}

\newcommand{\nv}{\textbf{n}}

\newcommand{\zrb}{\bm{0}}

\newcommand{\psib}{\bm{\psi}}

\newcommand{\Lmdb}{\bm{\Lambda}}
\newcommand{\Smb}{\bm{\Sm}}

\newcommand{\ds}{\,\mathrm{d}s}
\newcommand{\dt}{\,\mathrm{d}t}

\newcommand{\q}{\quad}
\newcommand{\qq}{\qquad}
\newcommand{\qqq}{\qquad\quad}

\newcommand{\qqqqq}{\qquad\qquad\quad}

\newtheorem{remark}{Remark}[section]

\begin{document}

\title{A Stabilized Time-Domain Combined Field Integral Equation Using the Quasi-Helmholtz Projectors}

\author{Van~Chien~Le,~\IEEEmembership{Member,~IEEE,}
        Pierrick~Cordel,
        Francesco~P.~Andriulli,~\IEEEmembership{Senior~Member,~IEEE,}
        and~Kristof~Cools% <-this % stops a space
\thanks{Manuscript received April 19, 2005; revised August 26, 2015. This work was supported by the European Research Council (ERC) under the European Union’s Horizon 2020 research and innovation programme (Grant agreement No. 101001847) and by the H2020-MSCA-ITN-EID project COMPETE GA No. 955476.}%
\thanks{Van Chien Le and Kristof Cools are with IDLab, Department of Information Technology at Ghent University - imec, 9000 Ghent, Belgium (e-mail: vanchien.le@ugent.be, kristof.cools@ugent.be).}% <-this % stops a space
\thanks{Pierrick Cordel and Francesco P. Andriulli are with the Department of Electronics and Telecommunications, Politecnico di Torino, 10129 Turin, Italy.}}% <-this % stops a space

\markboth{IEEE TRANSACTIONS ON ANTENNAS AND PROPAGATION}{LE \MakeLowercase{\textit{et al.}}: A stabilized TD-CFIE using the quasi-Helmholtz projectors}

\maketitle

\begin{abstract}
This paper introduces a time-domain combined field integral equation for electromagnetic scattering by a perfect electric conductor. The new equation is obtained by leveraging the quasi-Helmholtz projectors, which separate both the unknown and the source fields into solenoidal and irrotational components. These two components are then appropriately rescaled to cure the solution from a loss of accuracy occurring when the time step is large. Yukawa-type integral operators of a purely imaginary wave number are also used as a Calder\'{o}n preconditioner to eliminate the ill-conditioning of matrix systems. The stabilized time-domain electric and magnetic field integral equations are linearly combined in a Calder\'{o}n-like fashion, then temporally discretized using an appropriate pair of trial functions, resulting in a marching-on-in-time linear system. The novel formulation is immune to spurious resonances, dense discretization breakdown, large-time step breakdown and dc instabilities stemming from non-trivial kernels. Numerical results for both simply-connected and multiply-connected scatterers corroborate the theoretical analysis.
\end{abstract}


% \IEEEpeerreviewmaketitle

\section{Introduction}
\IEEEPARstart{T}{his} paper concerns the electromagnetic scattering of a transient plane wave by perfect electric conductors (PECs), which is commonly modeled by the time-domain electric field integral equation (TD-EFIE) and the time-domain magnetic field integral equation (TD-MFIE). One of the most prominent numerical methods to solve these equations is the marching-on-in-time (MOT) boundary element scheme. This scheme requires a discretization procedure in both space and time, resulting in matrix systems. Unfortunately, the TD-EFIE and TD-MFIE, upon discretization, suffer from several issues, which prevent scientists and engineers from developing a stable and efficient numerical scheme, as well as from obtaining accurate solution for all frequency ranges and different geometries.

First of all, the TD-EFIE systems become ill-conditioned when the spatial discretization is dense and/or the time step is large. These ill-conditioning are termed dense discretization breakdown and large-time step breakdown (low-frequency breakdown), respectively \cite{CAO+2009b}. Moreover, the TD-EFIE admits a non-trivial nullspace consisting of all constant-in-time solenoidal functions. This nullspace is the origin of direct current (dc) instabilities, which even lead the solution of the TD-EFIE to growing exponentially at late times \cite{ACO+2009,SBL2014}.

The TD-MFIE, in contrast, yields well-conditioned systems in both regimes and does not support any non-trivial kernel. However, the radar cross section computed from the MFIE's solution is less accurate than that computed from the EFIE's solution \cite{Peterson2008,GE2009}. Moreover, the MOT solution to the TD-MFIE on a toroidal surface contains slowly-oscillating non-decaying error components. These errors are termed near-dc instabilities and originated from the nullspaces of the static inner and outer MFIE operators \cite{CAO+2009a}.

In addition to the aforementioned issues, both the TD-EFIE and TD-MFIE suffer from resonant instabilities, which manifest themselves in harmonically oscillating error components superimposed on the MOT solutions and most prominently present at late times \cite{ACO+2009}. Besides that, when the time step $\Delta t$ is large, the scaling of solenoidal and irrotational components of the solution and the source fields differ by a factor of $\Delta t$. This means that, in practice, components of smaller order will be dominated or even suffer from a loss of significant digits, leading to completely wrong results in the far-field computation (see \cite{ZC2000,ZCC+2003,QC2010,BCA+2014} for more details). 

Most individual numerical issues and their combinations have been intensively studied and effectively resolved. For instance, the ill-conditioning of the TD-EFIE can be remedied using a Calder\'{o}n preconditioning strategy \cite{ACO+2009,CAO+2009b}, and dc instabilities can be tackled by applying the dot-trick formulation \cite{CAO+2009b} or the loop-tree decomposition \cite{WPC+2004,PW2005}. The inaccuracy of the TD-MFIE's solution has been mitigated by leveraging a mixed discretization scheme \cite{UBC+2017}. Resonant instabilities are usually eliminated by considering a time-domain combined field integral equation (TD-CFIE) \cite{SEA+2000,ACO+2009,BCB+2013}. 

The low-frequency breakdown and the loss of accuracy can be overcome by means of the loop-star/loop-tree decompositions, which are particularly effective with simply-connected scatterers \cite{PW2005,YJE2007,TXF2015}. The quasi-Helmholtz projectors have been introduced in \cite{ACB+2013} as an alternative to those decompositions. Numerous recent efforts leveraging the quasi-Helmholtz projectors to stabilize the low-frequency electromagnetic problems, both in the frequency-domain regime \cite{ACB+2013,MBC+2020,HEA+2023,HEA+2023b} and the time-domain regime \cite{DAC2020,BCA2015,BCA2015b}, have shown their effectiveness when applied to both simply-connected and multiply-connected geometries. In addition to the standard integral equations for the unknown surface electric current, some other approaches have also been proposed to deal with the low-frequency breakdown, such as the current-charge formulation \cite{TY2006}, the decoupled potential formulation \cite{VFG+2016}, and the Debye source approach \cite{EG2010,CB2013}.

Nonetheless, a formulation resolving all above mentioned issues is still missing from literature. This deficiency can be explained by the fact that PEC cavities support only discrete interior resonances (i.e., the lowest resonant frequency is far from the very-low-frequency band). Therefore, it does not make sense to consider extremely low frequencies and resonant frequencies in the same regime. However, the loss of accuracy occurs even at moderately low frequencies \cite{MBC+2020}. In this paper, we propose a stabilized TD-CFIE formulation for moderate to extremely low frequencies, which is immune to resonant instabilities, dense discretization breakdown and large-time step breakdown. For simply-connected geometries, the method is also free from dc instabilities. For multiply-connected geometries, near-dc instabilities are present, but correspond to exponentially decaying error components. Numerical experiments demonstrate they are not problematic in typical applications.

For this purpose, we extend the idea of the quasi-Helmholtz projected (qHP) frequency-domain CFIE introduced in \cite{MBC+2020}. Firstly, the solenoidal and irrotational components of the unknown and the incident fields are separated using the quasi-Helmholtz projectors and rescaled by appropriate factors, resulting in auxiliary well-balanced unknown and source fields. Secondly, the rescaled TD-EFIE and TD-MFIE are preconditioned by the corresponding stabilized Yukawa-type integral operators of a purely imaginary wave number. Then, their linearly combined formulation is temporally discretized using an appropriate pair of trial functions. The resulting linear system is solved by the MOT algorithm, and a post-processing procedure is performed to obtain the physical unknown.

The paper is organized as follows. The next section gives a brief introduction to the TD-EFIE, the TD-MFIE and Yukawa-type integral operators, together with their appropriate spatial discretization. This section also introduces the loop-star decomposition and the quasi-Helmholtz projectors associated with the discrete spatial function spaces. Section~\ref{sec:scaling} investigates the low-frequency (large-time step) behavior of the plane-wave incident fields, the induced current density and integral operators, as well as their stabilization using the quasi-Helmholtz projectors. Section~\ref{sec:qHP_TDCFIE} is the central part of this paper, which proposes a stabilized TD-CFIE and its proper temporal discretization. Some numerical results are presented in Section~\ref{sec:results} to corroborate our theoretical analysis. Finally, concluding remarks are provided in Section~\ref{sec:conclusion}.

\section{Background and notations}
\label{sec:background}

\subsection{Time-domain integral equations}

Let $\Gm$ be the surface of a PEC immersed in a homogeneous background medium with permittivity $\epsilon$ and permeability $\mu$. We assume that the surface $\Gm$ is closed, connected and orientable. The outward unit normal vector to $\Gm$ is denoted by $\nv$. An incident plane wave $(\eb^{in}, \hb^{in})$ induces a surface current density $\jb$ on $\Gm$, which satisfies the following standard TD-EFIE and TD-MFIE
\begin{align}
    & \,\, (\TT \jb)(\rb, t) = -\nv \times \eb^{in}(\rb, t), \label{eq:TD_EFIE} \\
    & \paren{\dfrac{1}{2} I + \KK}\jb(\rb, t) = \nv \times \hb^{in}(\rb, t), \label{eq:TD_MFIE}
\end{align}
for $\rb \in \Gm$ and $t > 0$. Here, $I$ stands for the identity operator and the time-domain integral operators $\TT$ and $\KK$ are defined by
\begin{align*}
    (\TT \jb)(\rb, t) & = (\TT^s \jb)(\rb, t) + (\TT^h \jb)(\rb, t), \\  
    (\TT^s \jb)(\rb, t) & = - \dfrac{\eta}{c} \, \nv \times \int_{\Gm} \dfrac{\pa_t \jb(\rb^\prime, \tau)}{4\pi R} \ds^\prime, \\
    (\TT^h \jb)(\rb, t) & = c \eta \, \nv \times \gradt_{\xxs} \int_{\Gm} \int_{-\infty}^{\tau}\dfrac{\divt_\Gm \, \jb(\rb^\prime, t^\prime)}{4\pi R} \dt^\prime \ds^\prime, \\
    (\KK\jb)(\rb, t) & = -\nv \times \curlt_{\xxs} \int_{\Gm} \dfrac{\jb(\rb^\prime, \tau)}{4\pi R} \ds^\prime,
\end{align*}
where $\eta = \sqrt{\mu/\epsilon}, c = 1/\sqrt{\mu\epsilon}$, $R = \abs{\rb - \rb^\prime}$ and $\tau = t - R/c$. In order to guarantee the uniqueness of a time-domain solution $\jb$, the causality is imposed, i.e., $\jb(\rb, t) = \zrb$ for all $t < 0$ and $\rb$ in a neighborhood of $\Gm$.

\subsection{Yukawa-type integral operators}

By a Fourier transform of the time-domain operators $\TT$ and $\KK$ at angular frequency $\omega > 0$, the corresponding frequency-domain operators read as
\begin{align*}
    (T_\kappa \jb)(\rb) & =  (T^s_{\kappa} \jb)(\rb) + (T^h_{\kappa} \jb)(\rb), \\
    (T^s_{\kappa} \jb)(\rb) & = -j \kappa \eta \, \nv \times \int_\Gm \dfrac{\exp(-j\kappa R)}{4\pi R} \, \jb(\rb^\prime) \ds^\prime, \\
    (T^h_{\kappa} \jb)(\rb) & = \dfrac{\eta}{j \kappa} \nv \times \gradt_{\xxs} \int_\Gm \dfrac{\exp(-j\kappa R)}{4\pi R} \, \divt_\Gm \, \jb(\rb^\prime) \ds^\prime, \\
    (K_{\kappa} \jb)(\rb) & = -\nv \times \curlt_{\xxs} \int_{\Gm} \dfrac{\exp(-j\kappa R)}{4\pi R} \, \jb(\rb^\prime) \ds^\prime,
\end{align*}
where $\kappa = \omega/c > 0$ is the wave number and $j = \sqrt{-1}$ is the imaginary unit. Based on this definition, in this paper, we consider the ``frequency-domain'' operators of purely imaginary wave number $-j\kappa$ with $\kappa > 0$, namely $T_{-j\kappa}$ and $K_{-j\kappa}$. These operators are usually called the Yukawa-type integral operators because they are solutions to a Yukawa-type equation \cite{SW2009}.

For the purpose of rendering the Yukawa-type integral operators $T_{-j\kappa}$ and $K_{-j\kappa}$ scaling identically to the time-domain operators $\TT$ and $\KK$ when the time step $\Delta t$ is large, the parameter $\kappa$ is chosen as $\kappa = (c\Delta t)^{-1}$ (see Sect.~\ref{sec:scaling} for more details on their large-time step scaling). The reader is also referred to \cite{LCA+2023} for the performance of the Yukawa-Calder\'{o}n TD-CFIE formulation corresponding to this value of $\kappa$. 

\begin{remark}
    In the large-time step regime (considered in this paper), the choice $\kappa = (c\Delta t)^{-1}$ (or $\kappa$ of the order of $(c\Delta t)^{-1}$) resolves the large-time step breakdown. In the small-time step regime, other problems occur, and more general choices for the regularizing operators could be appropriate \cite{EAG2014,CLM+2023}.
\end{remark}

\subsection{Spatial discretization}

\begin{figure}
    \centering
    \begin{tikzpicture}
      \draw[line width=0.4mm] (0, 0) -- (1, 1.73205) -- (2, 0) -- (1, -1.73205) -- (0, 0);
      \draw[line width=0.4mm, fill=gray!40] (2, 0) -- (3, 1.73205) -- (4, 0) -- (3, -1.73205) -- (2, 0);
      \draw[line width=0.4mm] (4, 0) -- (5, 1.73205) -- (6, 0) -- (5, -1.73205) -- (4, 0);
      \draw[line width=0.4mm] (1, 1.73205) -- (5, 1.73205);
      \draw[line width=0.4mm] (1, -1.73205) -- (5, -1.73205);
      \draw[line width=0.4mm] (0, 0) -- (6, 0);
      \draw[-latex, line width=0.6mm] (3, -0.6) -- (3, 0.6);
      \draw[fill=black] (3,1.73205) circle (2pt);
      \draw[fill=black] (3,-1.73205) circle (2pt);
      \draw[fill=black] (2,0) circle (2pt);
      \draw[fill=black] (4,0) circle (2pt);
      \node at (3.5, 0.2) {$e_n$};
      \node at (3, 2.08) {$\rb_n^-$};
      \node at (3, -2.05) {$\rb_n^+$};
      \node at (3, 1) {$c_n^-$};
      \node at (3, -0.95) {$c_n^+$};
    \end{tikzpicture}
    \caption{Rao-Wilton-Glisson (RWG) function $\fb_n(\rb)$ associated with the edge $e_n$. The function $\fb_n(\rb)$ is defined on two adjacent triangles $c^-_n$ and $c^+_n$ sharing the common edge $e_n$. The bold arrow indicates the direction of $\fb_n(\rb)$ across $e_n$.}
    \label{fig:RWG_function}
\end{figure}
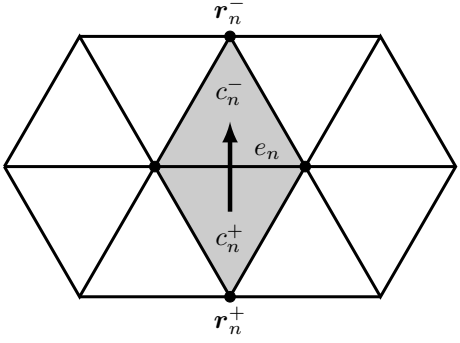

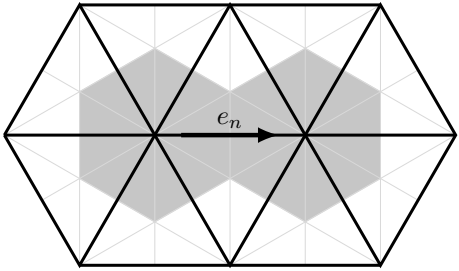
\begin{figure}
    \centering
    \begin{tikzpicture}
      \draw[line width=0.05mm, gray!30, fill=gray!45] (1, 0.57735) -- (1, -0.57735) -- (2, -1.1547) -- (3, -0.57735) -- (4, -1.1547) -- (5, -0.57735) -- (5, 0.57735) -- (4, 1.1547) -- (3, 0.57735) -- (2, 1.1547) -- (1, 0.57735) -- (1, 0.57735);
      \draw[line width=0.05mm, gray!30] (1, 1.73205) -- (1, -1.73205);
      \draw[line width=0.05mm, gray!30] (2, 1.73205) -- (2, -1.73205);
      \draw[line width=0.05mm, gray!30] (3, 1.73205) -- (3, -1.73205);
      \draw[line width=0.05mm, gray!30] (4, 1.73205) -- (4, -1.73205);
      \draw[line width=0.05mm, gray!30] (5, 1.73205) -- (5, -1.73205);
      \draw[line width=0.05mm, gray!30] (0, 0) -- (3, 1.73205);
      \draw[line width=0.05mm, gray!30] (6, 0) -- (3, 1.73205);
      \draw[line width=0.05mm, gray!30] (0, 0) -- (3, -1.73205);
      \draw[line width=0.05mm, gray!30] (6, 0) -- (3, -1.73205);
      \draw[line width=0.05mm, gray!30] (0.5, 0.866025) -- (5, -1.73205);
      \draw[line width=0.05mm, gray!30] (0.5, -0.866025) -- (5, 1.73205);
      \draw[line width=0.05mm, gray!30] (5.5, -0.866025) -- (1, 1.73205);
      \draw[line width=0.05mm, gray!30] (5.5, 0.866025) -- (1, -1.73205);
      \draw[line width=0.4mm] (0, 0) -- (1, 1.73205) -- (2, 0) -- (1, -1.73205) -- (0, 0);
      \draw[line width=0.4mm] (2, 0) -- (3, 1.73205) -- (4, 0) -- (3, -1.73205) -- (2, 0);
      \draw[line width=0.4mm] (4, 0) -- (5, 1.73205) -- (6, 0) -- (5, -1.73205) -- (4, 0);
      \draw[line width=0.4mm] (1, 1.73205) -- (5, 1.73205);
      \draw[line width=0.4mm] (1, -1.73205) -- (5, -1.73205);
      \draw[line width=0.4mm] (0, 0) -- (6, 0);
      \draw[-latex, line width=0.6mm] (2.35, 0) -- (3.65, 0);
      \node at (3, 0.2) {$e_n$};
    \end{tikzpicture}
    \caption{Buffa-Christiansen (BC) function $\gb_n(\rb)$ associated with the edge $e_n$. The function $\gb_n(\rb)$ is a linear combination of RWG functions defined on the barycentric refinement of the original mesh. The gray area indicates the support of $\gb_n(\rb)$, while the bold arrow indicates the direction of $\gb_n(\rb)$ along $e_n$.}
    \label{fig:BC_function}
\end{figure}

In order to numerically solve boundary integral equations, the surface $\Gm$ is partitioned into $N_C$ planar triangles with $N_E$ edges and $N_V$ vertices. On the triangular mesh of $\Gm$, we equip two boundary element spaces of divergence-conforming functions, namely Rao-Wilton-Glisson (RWG) functions $\fb_n(\rb)$ and Buffa-Christiansen (BC) functions $\gb_n(\rb)$, with $n = 1, 2, \ldots, N_E$. The RWG function $\fb_n(\rb)$ associated with the edge $e_n$ is defined on two adjacent triangles $c^+_n$ and $c^-_n$, which share the common edge $e_n$, as follows \cite{RWG1982}
\[
    \fb_n(\rb) = 
    \begin{cases}
        \dfrac{\rb - \rb^+_n}{2A_{c^+_n}} \qq & \text{for } \rb \in c^+_n, \\
        \dfrac{\rb^-_n - \rb}{2A_{c^-_n}} & \text{for } \rb \in c^-_n, \\
        0 & \text{otherwise},
    \end{cases}
\]
where $A_{c^+_n}$ and $A_{c^-_n}$ are the area of $c^+_n$ and $c^-_n$, respectively, while $\rb^+_n$ and $\rb^-_n$ respectively are the vertex of $c^+_n$ and $c^-_n$ that are opposite to $e_n$ (see Fig.~\ref{fig:RWG_function}). In other words, the function $\fb_n(\rb)$ represents a current flowing from $c^+_n$ to $c^-_n$. The BC function $\gb_n(\rb)$, which is associated with the edge $e_n$, is constructed as a linear combination of the RWG functions defined on the barycentric refinement of the original mesh, with the coefficients specified in \cite{BC2007} (see Fig.~\ref{fig:BC_function}).

The current density $\jb$ is spatially approximated as an expansion in the RWG basis functions
\[
    \jb(\rb, t) = \sum\limits_{n = 1}^{N_E} \jrm_n(t) \fb_n(\rb).
\]
The TD-EFIE \eqref{eq:TD_EFIE} is tested with the rotated RWG functions $\nv \times \fb_m(\rb)$, whereas the TD-MFIE \eqref{eq:TD_MFIE} is tested with the rotated BC functions $\nv \times \gb_m(\rb)$ \cite{CAD+2011,UBC+2017}, with $m = 1, 2, \ldots, N_E$, resulting in the following semi-discrete equations
\begin{align}
    & \,\,\, (\ZZ \jrm)(t) = -\erm(t), \label{eq:dis_TD_EFIE} \\
    & \paren{\dfrac{1}{2} \G + \MM} \jrm(t) = \hrm(t), \label{eq:dis_TD_MFIE}
\end{align}
where $\jrm(t) = \left[\jrm_1(t), \jrm_2(t), \ldots, \jrm_{N_E}(t) \right]^\transpose$, the Gram matrix
\[
    \left[\G\right]_{mn} = \inprod{\nv \times \gb_m, \fb_n},
\]
and the operators
\begin{align*}
    \ZZ & = \ZZ^s + \ZZ^h, \\
    \left[(\ZZ^s \jrm)(t) \right]_m & = \sum\limits_{n=1}^{N_E} \inprod{\nv \times \fb_m, \TT^s \paren{\jrm_n(t) \fb_n}}, \\
    \left[(\ZZ^h \jrm)(t) \right]_m & = \sum\limits_{n=1}^{N_E} \inprod{\nv \times \fb_m, \TT^h \paren{\jrm_n(t) \fb_n}}, \\
    \left[(\MM \jrm)(t) \right]_m & = \sum\limits_{n=1}^{N_E} \inprod{\nv \times \gb_m, \KK \paren{\jrm_n(t) \fb_n}},     
\end{align*}
with
\[
    \inprod{\fb, \gb} = \int_{\Gm} \fb(\rb) \cdot \gb(\rb) \ds.
\]
In addition, the semi-discrete right-hand sides are given by
\begin{align*}
    \left[\erm(t)\right]_m & = \inprod{\nv \times \fb_m, \nv \times \eb^{in}(t)}, \\
    \left[\hrm(t)\right]_m & = \inprod{\nv \times \gb_m, \nv \times \hb^{in}(t)}.
\end{align*}

Besides the time-domain integral operators, the Yukawa-type operators are also needed to be spatially discretized. In the next sections, the operators $T_{-j\kappa}$ and $K_{-j\kappa}$ are respectively composed with $\TT$ and $\KK$. To guarantee the stability of the composition, we discretize $T_{-j\kappa}$ using the BC functions, whereas $K_{-j\kappa}$ is discretized using the RWG trial functions and the rotated BC testing functions \cite{MBC+2020}. This scheme results in the following discretization matrices
\begin{align*}
    \Z & = \Z^s + \Z^h, \\
    \left[\Z^s\right]_{mn} & = \inprod{\nv \times \gb_m, T^s_{-j\kappa} \, \gb_n}, \\
    \left[\Z^h\right]_{mn} & = \inprod{\nv \times \gb_m, T^h_{-j\kappa} \, \gb_n}, \\
    \left[\M\right]_{mn} & = \inprod{\nv \times \gb_m, K_{-j\kappa} \fb_n}.
\end{align*} 

\subsection{Quasi-Helmholtz projectors}

This section is devoted to introducing the loop-star decomposition and the derived quasi-Helmholtz projectors for discrete spatial function spaces, which are essential in the stabilization of integral equations and operators throughout the paper. Every function $\jb$ in the solution space admits the decomposition \cite{Nedelec2001}
\[
    \jb = \curlt_{\Gm} \xi + \mathbf{\nabla}_{\Gm} \vphi + \psib, \qq \Delta_{\Gm} \psib = \zrb.
\]
Here, $\mathbf{\nabla}_{\Gm} \vphi$ is the irrotational component (star functions), while the local loop $\curlt_{\Gm} \xi$ and the global loop (harmonic functions) $\psib$ both represent solenoidal components. The space of harmonic functions has dimension $N_H := N_{holes} + 2 N_{handles}$, with $N_{holes}$ and $N_{handles}$, respectively, the number of holes and handles of $\Gm$. If the scatterer is simply-connected then $N_H = 0$.
In the RWG basis, local loops can be constructed via the vertex-edge connectivity matrix $\Lmdb \in \R^{N_E \times N_V}$ \cite{Vecchi1999,LLB2003,Andriulli2012}. Let us assume that the edge $e_m$ is oriented from $v^-_m$ to $v^+_m$, then
\[
    \left[\Lmdb\right]_{mn} :=
    \begin{cases}
        1 \q & \text{if the vertex } n \text{ is identical with } v^+_m, \\
        -1 & \text{if the vertex } n \text{ is identical with } v^-_m, \\
        0 & \text{otherwise}.
    \end{cases}
\]
On the other hand, RWG-based star functions can be constructed via the face-edge connectivity matrix $\Smb \in \R^{N_E \times N_C}$, defined by
\[
    \left[\Smb\right]_{mn} :=
    \begin{cases}
        1 \q & \text{if the triangle } n \text{ is identical with } c^+_m, \\
        -1 & \text{if the triangle } n \text{ is identical with } c^-_m, \\
        0 & \text{otherwise}.
    \end{cases}
\]
More specifically, the following quasi-Helmholtz decomposition holds for any RWG coefficient vector $\ib$
\[
    \ib = \Lmdb \lb + \Smb \ssb + \HHs \hb,
\]
where $\HHs \in \R^{N_E \times N_H}$. The vectors $\Lmdb \lb, \HHs \hb$ and $\Smb \ssb$ represent the RWG coefficients of local loops, global loops and star functions, respectively.

In the BC basis, given an arbitrary coefficient vector $\ovl{\ib}$, there exist three vectors $\ovl{\lb}, \ovl{\ssb}$ and $\ovl{\hb}$ such that \cite{SL2009,YJN2010,ACB+2013}
\[
    \ovl{\ib} = \Lmdb \ovl{\lb} + \Smb \ovl{\ssb} + \HHs \ovl{\hb}.
\]
In contrast to the RWG basis, $\Lmdb \ovl{\lb}$ and $\Smb \ovl{\ssb}$ are the BC coefficients of star functions and local loops, respectively.

Unfortunately, the construction of the standard quasi-Helmholtz decomposition is computationally expensive and unstable \cite{Andriulli2012}. The quasi-Helmholtz projectors were introduced in \cite{ACB+2013} as an alternative to cure the EFIE systems from low-frequency breakdown. More specifically, the quasi-Helmholtz projectors mapping RWG-based functions into the irrotational and solenoidal subspaces are respectively defined by
\begin{align*}
    \PPo^\Sm & = \Smb (\Smb^{\transpose} \Smb)^+ \Smb^\transpose, \\
    \PPo^{\Lambda H} & = I - \PPo^\Sm,
\end{align*}
where $(\Smb^{\transpose} \Smb)^+$ stands for the Moore–Penrose pseudoinverse matrix of $\Smb^{\transpose} \Smb$. Analogously, the quasi-Helmholtz projectors associated with the BC basis are given by
\begin{align*}
    \Pbb^\Lambda & = \Lmdb (\Lmdb^{\transpose} \Lmdb)^+ \Lmdb^\transpose, \\
    \Pbb^{\Sm H} & = I - \Pbb^\Lambda.
\end{align*}
Here, $\Pbb^\Lambda$ and $\Pbb^{\Sm H}$ are the projectors into the irrotational and solenoidal spaces, respectively. It is noteworthy that the quasi-Helmholtz projectors do not require an explicit identification of global loops $\HHs$. In addition, there is no distinction between local loop and global loop functions.

\section{Large-time step scaling and stabilization}
\label{sec:scaling}

In this section, the large-time step scaling of physical quantities and integral operators are analyzed. Those scaling are primarily derived from the corresponding low-frequency behavior of frequency-domain quantities and operators. Based on that analysis, appropriate stabilization are also proposed. Throughout this paper, all scaling matrices are with respect to the quasi-Helmholtz basis $\left[\Lmdb, \HHs, \Smb \right]^\transpose$.

\subsection{Plane-wave scattering}

When the time step $\Delta t$ is large, the discretized plane-wave incident fields $\erm$ and $\hrm$ scale as follows (see \cite{ZC2000,ZCC+2003,QC2010,BCA+2014,UBC+2017})
\[
    \erm = \OO
    \begin{pmatrix}
        \Delta t^{-1} \\
        \Delta t^{-1} \\
        1
    \end{pmatrix}, 
    \qq
    \hrm = \OO
    \begin{pmatrix}
        1 \\
        \Delta t^{-1} \\
        \Delta t^{-1}
    \end{pmatrix}.
\]
Please note that the scaling of $\erm$ and $\hrm$ are not identical because they correspond to different spatial discretization schemes. Next, the induced current density $\jrm$ has the following large-time step scaling \cite{ACB+2013}
\[
    \jrm = \OO
    \begin{pmatrix}
        1 \\
        1 \\
        \Delta t^{-1}
    \end{pmatrix}.
\]
This scaling implies that the irrotational components of the solution $\jrm$  (corresponding to star functions $\Smb$) are dominated by the solenoidal components when the time step $\Delta t$ is large. It leads to a loss of accuracy on the solution occurring even with moderately large time steps. In order to mitigate this loss, we introduce an auxiliary unknown \cite{BCA2015,BCA2015b}
\[
    \yrm(t) := \paren{\PPo^{\Lambda H} + \, \tilde{\pa}_t^{-1} \PPo^\Sm} \jrm(t),
\]
or
\begin{equation}
\label{eq:j}
    \jrm(t) = \paren{\PPo^{\Lambda H} + \, \tilde{\pa}_t \PPo^\Sm} \yrm(t),
\end{equation}
where $\tilde{\pa}_t^{-1} = \pa_t^{-1}/T_0$ and $\tilde{\pa}_t = T_0 \, \pa_t$, with $\pa_t^{-1}$ the temporal integral over $(-\infty, t)$, $T_0 = D/c$, and $D$ the diameter of the scatterer. The factor $T_0$ introduced in \cite{BCA2015b} aims to render the ``scaled'' integral $\tilde{\pa}_t^{-1}$ and the ``scaled'' derivative $\tilde{\pa}_t$ dimensionless (i.e., independent of the chosen units of measurement). Clearly, $\yrm = \OO(1, 1, 1)^\transpose$, whose components are well-balanced. Irrotational and solenoidal components of the incident electric field $\erm$ can also be balanced by using a multiplicative matrix as follows
\[
    \tilde{\erm}(t) := \paren{\widetilde{\Delta} t \PPo^{\Lambda H} + \PPo^\Sm} \erm(t).
\]
Here, the scaling factor $\widetilde{\Delta} t = \Delta t / T_0$ is also dimensionless. As a result, $\tilde{\erm} = \OO(1, 1, 1)^\transpose$. This multiplication only rescales the solenoidal components by a factor of $\widetilde{\Delta} t$, but it introduces neither differentiation nor integration in time. A stabilization procedure for $\hrm$ will be introduced in Sect.~\ref{subsec:MFIE}.

\subsection{Time-domain electric field integral equation}
\label{subsec:EFIE}

In the quasi-Helmholtz basis $\left[\Lmdb, \HHs, \Smb \right]^\transpose$, the large-time step scaling of the semi-discrete operator $\ZZ$ reads as \cite{ZC2000,ACB+2013}
\[
    \ZZ = \OO
    \begin{pmatrix}
        \Delta t^{-1} & \Delta t^{-1} & \Delta t^{-1} \\
        \Delta t^{-1} & \Delta t^{-1} & \Delta t^{-1} \\
        \Delta t^{-1} & \Delta t^{-1} & \Delta t
    \end{pmatrix}.
\]
Using the auxiliary unknown $\yrm$ and the balanced incident electric field $\tilde{\erm}$, we obtain the following stabilized TD-EFIE
\begin{equation}
    \label{eq:reg_TD_EFIE}
    (\widetilde{\ZZ} \yrm)(t) = - \tilde{\erm}(t),
\end{equation}
where the stabilized TD-EFIE operator
\begin{align*}
    \widetilde{\ZZ} 
    & := \paren{\widetilde{\Delta} t \PPo^{\Lambda H} + \PPo^\Sm} \ZZ \paren{\PPo^{\Lambda H} + \tilde{\pa}_t \PPo^\Sm} \\
    & \,\, = \paren{\widetilde{\Delta} t \PPo^{\Lambda H} + \PPo^\Sm} \ZZ^s \paren{\PPo^{\Lambda H} + \tilde{\pa}_t \PPo^\Sm} + \PPo^\Sm \tilde{\pa}_t \ZZ^h \PPo^\Sm.
\end{align*}
Here, we use the fact that $\PPo^{\Lambda H} \ZZ^h = \ZZ^h \PPo^{\Lambda H} = \zrb$. Equation \eqref{eq:reg_TD_EFIE} scales as 
\[    
    \begin{pmatrix}
        1 & 1 & \Delta t^{-1} \\
        1 & 1 & \Delta t^{-1} \\
        \Delta t^{-1} & \Delta t^{-1} & 1
    \end{pmatrix}
    \begin{pmatrix}
        1 \\
        1 \\
        1
    \end{pmatrix}
    = 
    \begin{pmatrix}
        1 \\
        1 \\
        1
    \end{pmatrix},
\]
which implies that \eqref{eq:reg_TD_EFIE} does not suffer from the large-time step breakdown and the loss of accuracy on the solution. However, like the TD-EFIE operator $\ZZ$, the operator $\widetilde{\ZZ}$ has a kernel comprising all constant-in-time solenoidal functions.

Analogously, the discretization matrix $\Z$ of the Yukawa-type EFIE operator $T_{-j\kappa}$, with $\kappa = (c\Delta t)^{-1}$, scales as
\[
    \Z = \OO
    \begin{pmatrix}
        \Delta t & \Delta t^{-1} & \Delta t^{-1} \\
        \Delta t^{-1} & \Delta t^{-1} & \Delta t^{-1} \\
        \Delta t^{-1} & \Delta t^{-1} & \Delta t^{-1}
    \end{pmatrix}.
\]
Please be aware that the role of $\Lmdb$ and $\Smb$ are exchanged from the RWG basis to the BC basis. Using the quasi-Helmholtz projectors $\Pbb^{\Sigma H}$ and $\Pbb^\Lambda$ of the BC functions, we can introduce the following stabilized Yukawa-type EFIE operator
\begin{align*}
    \widetilde{\Z} & := \paren{\widetilde{\Delta} t \, \Pbb^{\Sigma H} + \Pbb^\Lambda} \Z \paren{\Pbb^{\Sigma H} + \widetilde{\Delta} t^{-1} \Pbb^\Lambda} \\
    & = \paren{\widetilde{\Delta} t \, \Pbb^{\Sigma H} + \Pbb^\Lambda} \Z^s \paren{\Pbb^{\Sigma H} + \widetilde{\Delta} t^{-1} \Pbb^\Lambda} + \widetilde{\Delta} t^{-1} \Pbb^\Lambda \Z^h \Pbb^\Lambda,
\end{align*}
which has the large-time step scaling
\[    
    \widetilde{\Z} = \OO
    \begin{pmatrix}
        1 & \Delta t^{-1} & \Delta t^{-1} \\
        \Delta t^{-1} & 1 & 1 \\
        \Delta t^{-1} & 1 & 1
    \end{pmatrix}.
\]
Like $\ZZ^h$, we note that $\Pbb^{\Sigma H} \Z^h = \Z^h \Pbb^{\Sigma H} = \zrb$. Finally, the matrix $\widetilde{\Z}$ can be used as a Calder\'{o}n preconditioner of $\widetilde{\ZZ}$, which results in the following well-conditioned, large-time step well-balanced TD-EFIE
\begin{equation}
    \label{eq:qhp_CP_TDEFIE}
    \widetilde{\Z} \, \G^{-\transpose} \widetilde{\ZZ} \, \yrm(t) = - \widetilde{\Z} \, \G^{-\transpose} \tilde{\erm}(t).
\end{equation}
It is noteworthy that the preconditioner $\widetilde{\Z}$ is independent of time, thereby saving the cost of additional temporal discretization procedures, such as the scheme proposed in \cite{BCA2015b}.

\subsection{Time-domain magnetic field integral equation}
\label{subsec:MFIE}

In this section, we consider the following symmetrized TD-MFIE, which is inspired by the Calder\'{o}n-like frequency-domain MFIE formulation in \cite{MBC+2020,LC2023}
\begin{equation}
    \label{eq:sym_TD_MFIE}
    \X \, \G^{-1} \XX \, \jrm(t) = \X \, \G^{-1} \hrm(t),
\end{equation}
where 
\[
    \X := \paren{\dfrac{1}{2}\G - \M}, \q \text{and} \qq \XX := \paren{\frac{1}{2}\G + \MM}.
\]
The operator $\X \, \G^{-1} \XX$ is ``symmetric'' in the sense that $\X$ and $\XX$ are corresponding to the inner and outer MFIE operators, respectively. 
Since the operators $I$, $\KK$ and $K_{-j\kappa}$ are discretized using the same spatial discretization scheme, the semi-discrete operator $\XX$ and the matrix $\X$ share the following large-time step scaling (see \cite{ZCC+2003,BCA+2014,UBC+2017})
\[
    \XX \sim \X = \OO\begin{pmatrix}
         1 & 1 & 1 \\
        \Delta t^{-2} & 1 & 1 \\
        \Delta t^{-2} & \Delta t^{-2} & 1
    \end{pmatrix}.
\]
In addition, the algebraic structure of the inverse Gram matrix $\G^{-1}$ is as follows \cite{BCA2015b}
\[
    \G^{-1} =
    \begin{pmatrix}
        \square & \square & \square \\
        0 & \square & \square \\
        0 & 0 & \square
    \end{pmatrix}.
\]
Here, the square symbols represent non-zero blocks. Based on those behaviors, the large-time step scaling of the operator $\X \, \G^{-1} \XX$, the unknown $\jrm$ and the right-hand side $\X \, \G^{-1} \hrm$ in \eqref{eq:sym_TD_MFIE} respectively read as
\[
    \begin{pmatrix}
        1 & 1 & 1 \\
        \Delta t^{-2} & 1 & 1  \\
        \Delta t^{-2} & \Delta t^{-2} & 1
    \end{pmatrix}
    \begin{pmatrix}
        1 \\
        1 \\
        \Delta t^{-1}
    \end{pmatrix}
    \ \neq \
    \begin{pmatrix}
        1 \\
        \Delta t^{-1} \\
        \Delta t^{-1}
    \end{pmatrix}.  
\]
The symbol $\neq$ indicates that the scaling of the left-hand side and the right-hand side of \eqref{eq:sym_TD_MFIE} are not identical. Obviously, at least one of them is incorrect. This discrepancy can be eliminated by taking into account the following property of the static MFIE operators
\begin{equation}
    \label{eq:property}
    \Pbb^{\Sm H} \paren{\dfrac{1}{2}\G - \M_0} \G^{-1} \paren{\dfrac{1}{2}\G + \M_0} \PPo^{\Lambda H} = \zrb,
\end{equation}
where $\M_0$ is the discretization matrix of the static operator $K_0$. We refer the reader to \cite{MBC+2020} for a rigorous proof of \eqref{eq:property}. Besides that, for sufficiently large time step $\Delta t \gg D/c$, the following estimate can be deduced from the Maclaurin series of the exponential function
\begin{align*}
    (K_{-j\kappa} - K_0) \jb(\rb) 
    & = \nv \times \curlt_{\xxs} \int_{\Gm} \dfrac{1 - \exp(-\kappa R)}{4\pi R} \, \jb(\rb^\prime) \ds^\prime \\
    & \approx \nv \times \curlt_{\xxs} \int_{\Gm} \paren{\dfrac{\kappa}{4\pi} - \dfrac{ \kappa^2 R}{8\pi}} \jb(\rb^\prime) \ds^\prime \\
    & = - \nv \times \curlt_{\xxs} \int_{\Gm} \dfrac{ \kappa^2 R}{8\pi} \jb(\rb^\prime) \ds^\prime \\
    & = \OO(\kappa^2) = \OO(\Delta t^{-2}).
\end{align*}
Please keep in mind that $\kappa = (c\Delta t)^{-1}$. Moreover, it follows from the Taylor expansion of $\jb$ that
\begin{align*}
    (\KK - K_0) \jb(t) 
    & = \nv \times \curlt_{\xxs} \int_{\Gm} \dfrac{\jb(t) - \jb\paren{t - R/c}}{4\pi R} \ds^\prime \\
    & \approx \nv \times \curlt_{\xxs} \int_{\Gm} \paren{\dfrac{\pa_t\jb(t)}{4\pi c} - R\dfrac{\pa^2_t\jb(t)}{8\pi c^2}} \ds^\prime \\
    & = - \nv \times \curlt_{\xxs} \int_{\Gm} \dfrac{R}{8\pi c^2} \pa^2_t\jb(t) \ds^\prime \\
    & = \OO(\Delta t^{-2}).
\end{align*}
Here, the space variables $\rb$ and $\rb^\prime$ are omitted for the sake of readability. Now, we can rearrange the ``loop-loop'' components of the operator $\X \, \G^{-1} \XX$ as follows
\begin{align*}
    & \Pbb^{\Sm H} \paren{\dfrac{1}{2}\G - \M} \G^{-1} \paren{\dfrac{1}{2}\G + \MM} \PPo^{\Lambda H} \\
    & \q = \Pbb^{\Sm H} \paren{\M_0 - \M} \G^{-1} \paren{\MM - \M_0} \PPo^{\Lambda H} \\
    & \qq \, + \Pbb^{\Sm H} \paren{\dfrac{1}{2}\G - \M_0} \G^{-1} \paren{\MM - \M_0} \PPo^{\Lambda H} \\
    & \qq \, + \Pbb^{\Sm H} \paren{\M_0 - \M} \G^{-1} \paren{\dfrac{1}{2}\G + \M_0} \PPo^{\Lambda H} \\
    & \qq \, + \Pbb^{\Sm H} \paren{\dfrac{1}{2}\G - \M_0} \G^{-1} \paren{\dfrac{1}{2}\G + \M_0} \PPo^{\Lambda H} \\
    & \q = \OO(\Delta t^{-4}) + \OO(\Delta t^{-2}) + \OO(\Delta t^{-2}) + \zrb = \OO(\Delta t^{-2}).
\end{align*}
Consequently, we arrive at the following correct scaling of the operator $\X \, \G^{-1} \XX$
\[
    \X \, \G^{-1} \XX = \OO
    \begin{pmatrix}
        1 & 1 & 1 \\
        \Delta t^{-2} & \Delta t^{-2} & 1  \\
        \Delta t^{-2} & \Delta t^{-2} & 1
    \end{pmatrix}.
\]
Next, by using appropriate stabilization for the unknown $\jrm$ and the right-hand side $\widetilde{\X} \, \G^{-1} \hrm$, we end up with the following stabilized version of the symmetrized TD-MFIE \eqref{eq:sym_TD_MFIE}
\begin{equation}
    \label{eq:qHP_sym_TDMFIE}
    \widetilde{\X} \, \G^{-1} \widetilde{\XX} \, \yrm(t) = \widetilde{\X} \, \G^{-1} \hrm(t),
\end{equation}
where
\begin{align*}
    \widetilde{\X} & := \paren{\widetilde{\Delta} t \, \Pbb^{\Sigma H} + \Pbb^\Lambda} \X, \\
    \widetilde{\XX} & := \XX \paren{\PPo^{\Lambda H} + \tilde{\pa}_t \PPo^\Sm}.
\end{align*}
Equation \eqref{eq:qHP_sym_TDMFIE} properly scales as
\[
    \begin{pmatrix}
        1 & 1 & \Delta t^{-1} \\
        \Delta t^{-1} & \Delta t^{-1} & 1 \\
        \Delta t^{-1} & \Delta t^{-1} & 1
    \end{pmatrix}
    \begin{pmatrix}
        1 \\
        1 \\
        1
    \end{pmatrix}
    = 
    \begin{pmatrix}
        1 \\
        1 \\
        1
    \end{pmatrix},
\]
which implies that it is well-balanced when the time step $\Delta t$ is large. Because of the introduction of the auxiliary unknown $\yrm$, the stabilized TD-MFIE \eqref{eq:qHP_sym_TDMFIE} supports constant-in-time irrotational regime solutions.

\begin{remark}
    Even though the low-frequency (large-time step) and dense-mesh regimes share some spectral properties, dense-mesh behaviors of integral operators are more complicated than their large-time step scaling. In this paper, we leverage the Calder\'{o}n preconditioning technique to eliminate dense-mesh breakdown. This technique is commonly used and thoroughly analyzed in the literature. For the conciseness of the paper, we would then omit the dense-mesh analysis.
\end{remark}

\section{A stabilized TD-CFIE}
\label{sec:qHP_TDCFIE}

We are in the position to propose a stabilized TD-CFIE by linearly combining the TD-EFIE \eqref{eq:qhp_CP_TDEFIE} and the TD-MFIE \eqref{eq:qHP_sym_TDMFIE}. These two equations are compatible because they both map an RWG coefficient vector into the dual of the space spanned by the rotated BC functions. In particular, we consider the following combined formulation
\begin{equation}
    \label{eq:qHP_TDCFIE}
    \paren{\widetilde{\Z} \, \G^{-\transpose} \widetilde{\ZZ} + \eta^2 \, \widetilde{\X} \, \G^{-1} \widetilde{\XX}} \yrm = - \widetilde{\Z} \, \G^{-\transpose} \tilde{\erm} + \eta^2 \, \widetilde{\X} \, \G^{-1} \hrm.
\end{equation}
Obviously, \eqref{eq:qHP_TDCFIE} is large-time step well-balanced, which scales as follows
\[
    \begin{pmatrix}
        1 & 1 & \Delta t^{-1} \\
        1 & 1 & 1 \\
        1 & 1 & 1
    \end{pmatrix}
    \begin{pmatrix}
        1 \\
        1 \\
        1
    \end{pmatrix}
    = 
    \begin{pmatrix}
        1 \\
        1 \\
        1
    \end{pmatrix}.
\]

Moreover, it is corroborated by numerical experiments that the qHP TD-CFIE \eqref{eq:qHP_TDCFIE} does not suffer from harmonically-oscillating non-decaying errors stemming from interior resonances, and that \eqref{eq:qHP_TDCFIE} is well-conditioned when the spatial mesh is fine and the time step is large. 

It is noteworthy that the stabilized TD-CFIE \eqref{eq:qHP_TDCFIE} does not support any non-trivial kernel, since the kernel of the stablized TD-EFIE and TD-MFIE operators are disjoint (the former comprises constant-in-time solenoidal functions, whereas the latter comprises constant-in-time irrotational functions). It implies that \eqref{eq:qHP_TDCFIE} is immune to dc instabilities rooted in non-trivial nullspaces. Nonetheless, the polynomial eigenvalue analysis performed in Section~\ref{sec:results} reveals that for multiply-connected geometries, \eqref{eq:qHP_TDCFIE} suffers from near-dc instabilities, which are corresponding to exponentially-decaying error components.

In a relevant frequency-domain context, a direct formulation called Yukawa-Calder\'{o}n CFIE (a symmetrized combined formulation) has been introduced in \cite{MBC+2020}, where its dense grid stability, low-frequency stability, and unique solvability have been demonstrated. Its indirect counterpart has been studied in \cite{LC2023}, where the well-posedness on merely Lipschitz geometries, a coercivity estimate, and an alternative stable discretization have been put forward.

Clearly, the choice of the coupling parameter has a major impact on the spectral properties of \eqref{eq:qHP_TDCFIE}. In this paper, the optimal value $\eta^2$ is chosen, which minimizes the condition number of the system \eqref{eq:qHP_TDCFIE} for the electromagnetic scattering by the unit sphere \cite{LC2023,Kress1985}. Since the equation \eqref{eq:qHP_TDCFIE} is still time-dependent, a temporal discretization scheme is needed.

\subsection{Temporal discretization}

\begin{figure}
    \centering
    \begin{tikzpicture}
        \draw[-latex,line width=0.2mm] (-2.6,-0.3)--(3.3,-0.3) node[below=1mm]{$t$}; 
        \draw[-latex,line width=0.2mm] (-2.3,-0.6)--(-2.3,1.7) node[left]{$q_0(t)$};
        \draw (-1.2,-0.4)--(-1.2,-0.3) node[below=0.6mm]{$-\Delta t$};
        \draw (0,-0.4)--(0,-0.3) node[below=0.8mm]{$0$};
        \draw (1.2,-0.4)--(1.2,-0.3) node[below=0.6mm]{$\Delta t$};
        \draw (2.4,-0.4)--(2.4,-0.3) node[below=0.6mm]{$2\Delta t$};
        \draw (-2.3,0)--(-2.4,0) node[left]{$0$};
        \draw (-2.3,1.2)--(-2.4,1.2) node[left]{$1$};
        \draw[line width=0.45mm,black] (-1.9,0) -- (-1.2,0);
        \draw[domain=-1.2:0,smooth,variable=\x,line width = 0.45mm,black] plot ({\x},{0.6*(\x/1.2+1)^2});
        \draw[domain=0:1.2,smooth,variable=\x,line width = 0.45mm,black] plot ({\x},{0.6-\x^2/1.2+\x});
        \draw[domain=1.2:2.4,smooth,variable=\x,line width = 0.45mm,black] plot ({\x},{0.6*(\x/1.2-2)^2});
        \draw[line width=0.45mm,black] (2.4,0) -- (3,0);
    \end{tikzpicture}
    \caption{The continuously differentiable piecewise-quadratic function $q_0$.}
    \label{fig:quad_function}
\end{figure}
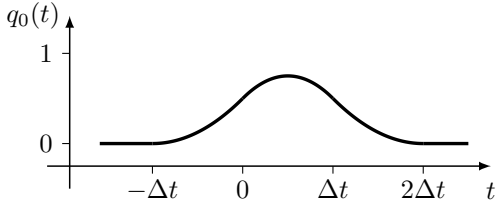

Let $N_T$ be the number of time steps that we want to compute the solution. In order to discretize the qHP TD-CFIE \eqref{eq:qHP_TDCFIE} in time, we introduce two set of temporal shifted basis functions
\[
    h_i(t) = h_0(t - i\Delta t), \qqq q_i(t) = q_0(t - i\Delta t),
\]
with $i = 1, 2, \ldots, N_T$. The continuous piecewise-linear (hat) function $h_0$ and the continuously differentiable piecewise-quadratic function $q_0$ are defined by (see Fig.~\ref{fig:hat_function} and Fig.~\ref{fig:quad_function})
\[
    h_0(t) = 
    \begin{cases}
        1 + \frac{t}{\Delta t} \qq & -\Delta t \le t < 0, \\
        1 - \frac{t}{\Delta t} &  0 \le t \le \Delta t, \\
        0 & \text{otherwise},
    \end{cases}
\]
and
\[
    q_0(t) = 
    \begin{cases}
        \frac{1}{2} \paren{\frac{t}{\Delta t} + 1}^2 \qq & -\Delta t \le t < 0, \\
        \frac{1}{2} - \frac{t^2}{\Delta t^2} + \frac{t}{\Delta t} &  0 \le t < \Delta t, \\
        \frac{1}{2} \paren{\frac{t}{\Delta t} - 2}^2 \qq & \Delta t \le t \le 2\Delta t, \\
        0 & \text{otherwise}.
    \end{cases}
\]

\begin{figure}
    \centering
    \begin{tikzpicture}
        \draw[-latex,line width=0.2mm] (-2.6,-0.3)--(3.3,-0.3) node[below=1mm]{$t$}; 
        \draw[-latex,line width=0.2mm] (-2.3,-0.6)--(-2.3,1.7) node[left]{$h_0(t)$};
        \draw (-1.2,-0.4)--(-1.2,-0.3) node[below=0.6mm]{$-\Delta t$};
        \draw (0,-0.4)--(0,-0.3) node[below=0.8mm]{$0$};
        \draw (1.2,-0.4)--(1.2,-0.3) node[below=0.6mm]{$\Delta t$};
        \draw (2.4,-0.4)--(2.4,-0.3) node[below=0.6mm]{$2\Delta t$};
        \draw (-2.3,0)--(-2.4,0) node[left]{$0$};
        \draw (-2.3,1.2)--(-2.4,1.2) node[left]{$1$};
        \draw[line width=0.45mm,black] (-2.1,0) -- (-1.2,0) -- (0,1.2) -- (1.2,0) -- (2.8,0);
    \end{tikzpicture}
    \caption{The continuous piecewise-linear function $h_0$.}
    \label{fig:hat_function}
\end{figure}
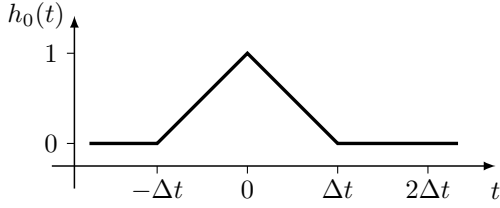

By introducing the auxiliary unknown $\yrm$ in \eqref{eq:j}, we introduce an additional differentiation in time to the irrotational components. It implies that the irrotational components of $\yrm$ must be discretized by temporal basis functions of one order higher than those used for the solenoidal components. More specifically, we demonstrate here the expansion of $\yrm$ in the lowest-order pair of temporal basis functions 
\begin{equation}
    \label{eq:y_expansion}
    \yrm(t) = \sum\limits_{i=1}^{N_T} \paren{h_i(t) \PPo^{\Lambda H} + \, q_i(t) \PPo^\Sm} \yys_i.
\end{equation}
Equation \eqref{eq:qHP_TDCFIE} is then tested with the Dirac delta distributions $\delta_k(t) = \delta(t - k\Delta t)$, with $k = 1, 2, \ldots, N_T$, resulting in the following lower triangular block matrix system
\begin{equation}
\label{eq:dis_TD_CFIE}
    \begin{pmatrix}
        \T_0 & & & \\
        \T_1 & \T_0 & & \\
        \vdots & \vdots & \ddots & \\
        \T_{N_T-1} & \T_{N_T-2} & \ldots & \T_0
    \end{pmatrix}   
    \begin{pmatrix}
        \yys_1 \\
        \yys_2 \\
        \vdots \\
        \, \yys_{N_T}
    \end{pmatrix} 
    =
    \begin{pmatrix}
        \rrs_1 \\
        \rrs_2 \\
        \vdots \\
        \rrs_{N_T}
    \end{pmatrix}.
\end{equation}
In \eqref{eq:dis_TD_CFIE}, the square matrix $\T_i$, with $i = 0, 1, \ldots, N_T - 1$, is defined by
\[
    \T_i := \widetilde{\Z} \, \G^{-\transpose} \widetilde{\ZZs}_i + \eta^2 \, \widetilde{\X} \, \G^{-1} \widetilde{\XXs}_i,
\]
where
\begin{align*}
     \widetilde{\ZZs}_i & := \int_{\R} \delta(t) \, \widetilde{\ZZ} \paren{h_i(t) \PPo^{\Lambda H} + \, q_i(t) \PPo^\Sm} \dt, \\
     \widetilde{\XXs}_i & := \int_{\R} \delta(t) \, \widetilde{\XX} \paren{h_i(t) \PPo^{\Lambda H} + \, q_i(t) \PPo^\Sm} \dt.
\end{align*}
Please note that the operators $\widetilde{\ZZ}$ and $\widetilde{\XX}$ do not contain any temporal integration. The right-hand side vector $\rrs_k$, with $k = 1, 2, \ldots, N_T$, is determined by
\[
    \rrs_k := - \widetilde{\Z} \, \G^{-\transpose} \int_\R \delta_k(t) \, \tilde{\erm}(t) \dt + \eta^2 \, \widetilde{\X} \, \G^{-1} \int_\R \delta_k(t) \, \hrm(t) \dt.
\]

\subsection{Marching-on-in-time algorithm}

The block matrix system \eqref{eq:dis_TD_CFIE} can be efficiently solved by the MOT algorithm
\[
    \yys_i = \T_0^{-1} \paren{\rrs_i - \sum\limits_{k=1}^{i-1} \T_k \yys_{i-k}}.
\]
When the coefficient vectors $\yys_i$ of the auxiliary unknown $\yrm$ are computed, we can use the definition \eqref{eq:j} and the expansion \eqref{eq:y_expansion} to get back the physical current density $\jrm(t)$ as follows
\begin{align*}
    \jrm(t) 
    & = \paren{\PPo^{\Lambda H} + \, \tilde{\pa}_t \PPo^\Sm} \yrm(t) \\
    & = \sum\limits_{i=1}^{N_T} \paren{h_i(t) \PPo^{\Lambda H} + \, \tilde{\pa}_t q_i(t) \PPo^\Sm} \yys_i \\
    & = \sum\limits_{i=1}^{N_T - 1} \paren{h_i(t) \PPo^{\Lambda H} + \dfrac{h_i(t) - h_{i+1}(t)}{\widetilde{\Delta} t} \PPo^\Sm} \yys_i \\
    & \qqqqq + \paren{\PPo^{\Lambda H} + \dfrac{1}{\widetilde{\Delta} t} \PPo^\Sm} \yys_{N_T} h_{N_T}(t) \\
    & = \sum\limits_{i=1}^{N_T} \paren{\PPo^{\Lambda H} \yys_i + \dfrac{1}{\widetilde{\Delta} t} \PPo^\Sm \paren{\yys_i -\yys_{i-1}}}  h_i(t).
\end{align*}
Here, we use the convention $\yys_0 = \zrb$. Denoting by $\jjs_i$ the expansion coefficient of $\jrm(t)$ corresponding to the temporal basis function $h_i(t)$, with $i = 1, 2, \ldots, N_T$, we arrive at
\begin{equation}
    \label{eq:post_processing}
    \jjs_i = \PPo^{\Lambda H} \yys_i + \dfrac{1}{\widetilde{\Delta} t} \PPo^\Sm \paren{\yys_i - \yys_{i-1}}.
\end{equation}

\begin{remark}
    The shifted piecewise-polynomial Lagrange interpolators of higher orders (also denoted by $h_i$) can be used to discretize the solenoidal components of $\yrm$ (see \cite{MMR1997} and \cite{BCB+2013} for some examples). For compatibility, the temporal basis functions $q_i$ for the irrotational components must be chosen such that
    \[
        q_i^\prime(t) = \dfrac{h_i(t) - h_{i+1}(t)}{\Delta t},
    \]  
    for all indices $i = 1, 2, \ldots, N_T-1$.
\end{remark}

\subsection{Complexity analysis}
The proposed qHP TD-CFIE involves rescaling and preconditioning procedures, necessitating additional matrix/vector product calculations compared to the standard TD-EFIE and TD-MFIE. Specifically, the rescaling involves the quasi-Helmholtz projectors, which can be constructed with linear complexity $\OO(N_E)$ \cite{ACB+2013}. Multiplication with these projectors also demands $\OO(N_E)$ operations. The preconditioning technique introduces the Gram matrix $\G$ and the Galerkin discretization matrices $\Z$ and $\M$ of the Yukawa-type operators $T_{-j\kappa}$ and $K_{-j\kappa}$. Whereas the former is sparse and well-conditioned, the latter are dense when the time step $\Delta t$ is large. Nevertheless, multiplication with them all have only a linear complexity $\OO(N_E)$, by means of acceleration techniques for quasi-static kernels \cite{Bebendorf2000}. Finally, the proposed MOT scheme can be coupled with fast time-domain algorithms, such as the plane-wave time-domain (PWTD) method \cite{ESM1999}, the time-domain adaptive integral method \cite{YJM2004}, the nonuniform grid time domain algorithm \cite{BLM2006}, and the accelerated Cartesian expansions \cite{VS2007}, thereby reducing the computational cost even to scaling quasi-linearly with the number of spatial unknowns $N_E$ and linearly with the number of temporal unknowns $N_T$. For instance, the total complexity of the proposed numerical scheme is $\OO(N_T N_E \log^2 N_E)$ using the PWTD method.

\section{Numerical results}
\label{sec:results}

\begin{figure}[htbp]
\centerline{
\includegraphics[width=0.3\linewidth]{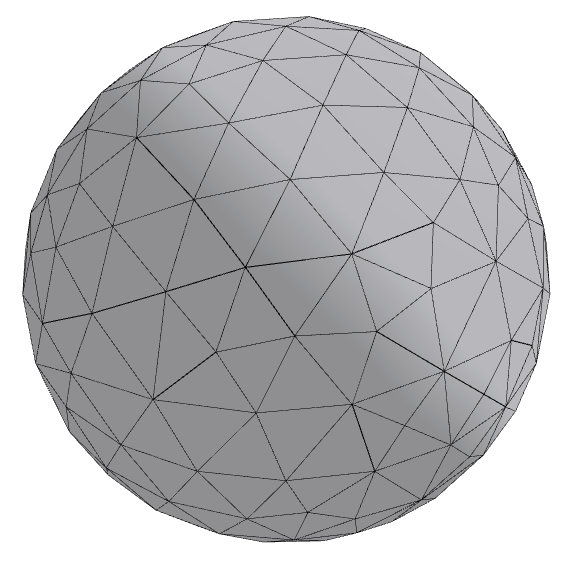} \qq
\includegraphics[width=0.43\linewidth]{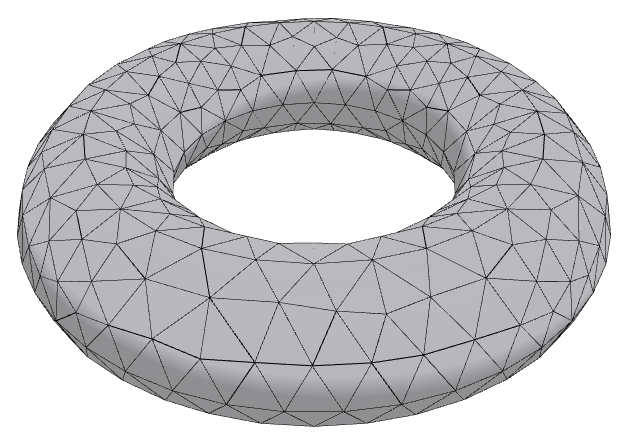}}
\caption{Triangular mesh of two surface geometries used in numerical experiments. \textit{Left:} a sphere of radius $1\mathrm{m}$. \textit{Right:} a torus of two radii $3\mathrm{m}$ and $1\mathrm{m}$.}
\label{fig:geometries}
\end{figure}

This section presents some numerical results for the scattering by PECs of two different geometries (see Fig.~\ref{fig:geometries}).
The excitation is provided by a Gaussian-in-time plane-wave electric field of the form
\[
    \eb^{in}(\rb, t) = \dfrac{4A}{w \sqrt{\pi}} \pb \exp \paren{-\paren{\dfrac{4}{w}\paren{c(t - t_0) - \kb \cdot \rb}}^2},
\]
with amplitude $A = 1\mathrm{V}$, polarization $\pb = \mathbf{1}_{x}$ and direction $\kb = \mathbf{1}_z$. The width $w$ and the time of arrival $t_0$ are specified on each experiment. The incident magnetic field $\hb^{in}$ is then determined by
\[
    \hb^{in}(\rb, t) = \eta^{-1} \kb \times \eb^{in}(\rb, t).
\]

In order to examine the performance of the qHP TD-CFIE \eqref{eq:qHP_TDCFIE}, its numerical results are compared with those of the following formulations
\begin{itemize}
    \item the standard TD-EFIE \eqref{eq:dis_TD_EFIE};
    \item the mixed discretized TD-MFIE \eqref{eq:dis_TD_MFIE};
    \item the mixed TD-CFIE proposed in \cite{BCB+2013}, with coupling parameter $\alpha = 0.5$.
\end{itemize}
These three equations are all temporally expanded in the set of hat functions $h_i(t)$, and then tested with the Dirac delta distributions $\delta_k(t)$. 

\subsection{Sphere}

We start with the scattering by a simple PEC scatterer: a sphere of radius $1\mathrm{m}$ (see Fig.~\ref{fig:geometries}, left). It is approximated by 476 triangles with 240 vertices and 714 edges. 

\begin{figure}[htbp]
\includegraphics[width=0.98\linewidth]{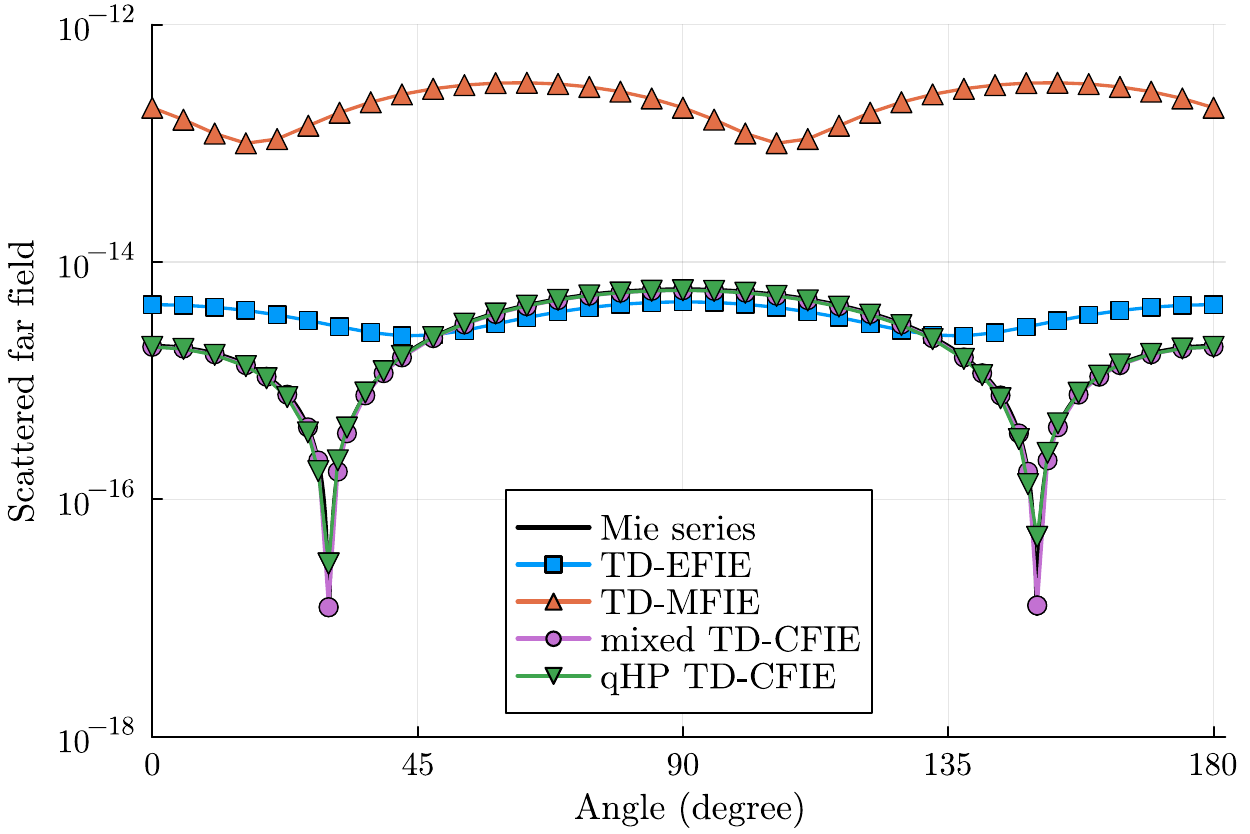}
\caption{Scattered far fields computed using the Mie series and different numerical schemes. Numerical solutions are computed with time step $\Delta t = 0.333 \mathrm{ms}$, then they are Fourier transformed at the frequency $3 \mathrm{Hz}$. The TD-EFIE's and TD-MFIE's solutions are inaccurate due to a loss of significant digits on the irrotational components at large time steps, whereas the solution to the qHP TD-CFIE is accurate.}
\label{fig:farfield}
\end{figure}

Firstly, the exemption from the loss of accuracy at large time steps is verified. The Gaussian-in-time plane wave with $w = 8\cdot 10^6\mathrm{m}$ and $t_0 = 80 \mathrm{ms}$ is used as the excitation. The unknown surface current $\jb(\rb, t)$ is numerically computed with time step $\Delta t = 0.333 \mathrm{ms}$, then it is Fourier transformed. Afterwards, the scattered far field computed from the transformed surface current at the frequency $3\mathrm{Hz}$ is compared with the reference one obtained by the Mie series. Fig.~\ref{fig:farfield} shows that when the time step $\Delta t$ is large, the TD-EFIE's and TD-MFIE's solutions exhibit a loss of accuracy, resulting in inaccurate far fields. In contrast, the qHP TD-CFIE does not suffer from that loss.

\begin{figure*}[!t]
\centering
\subfloat[]{\includegraphics[width=0.47\linewidth]{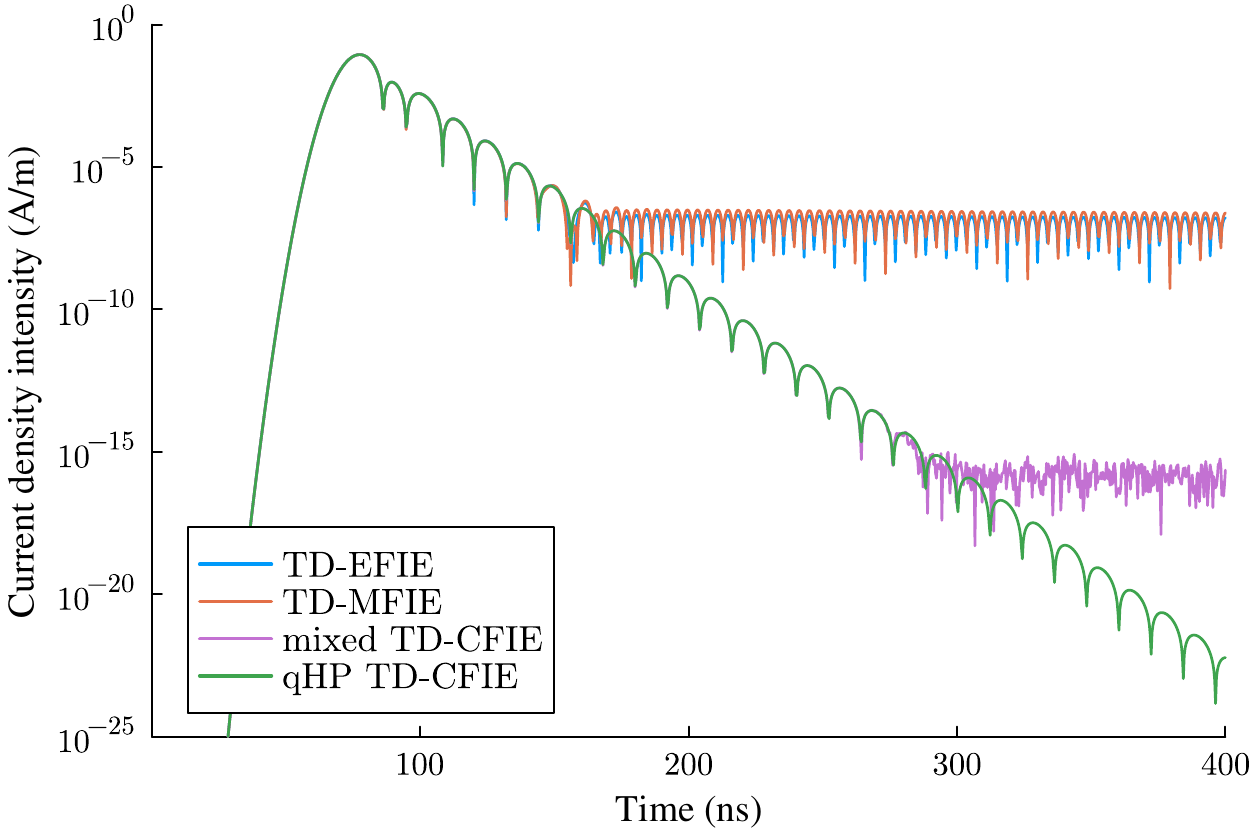}%
\label{fig:current_sphere}}
\hfil \hfil
\subfloat[]{\includegraphics[width=0.47\linewidth]{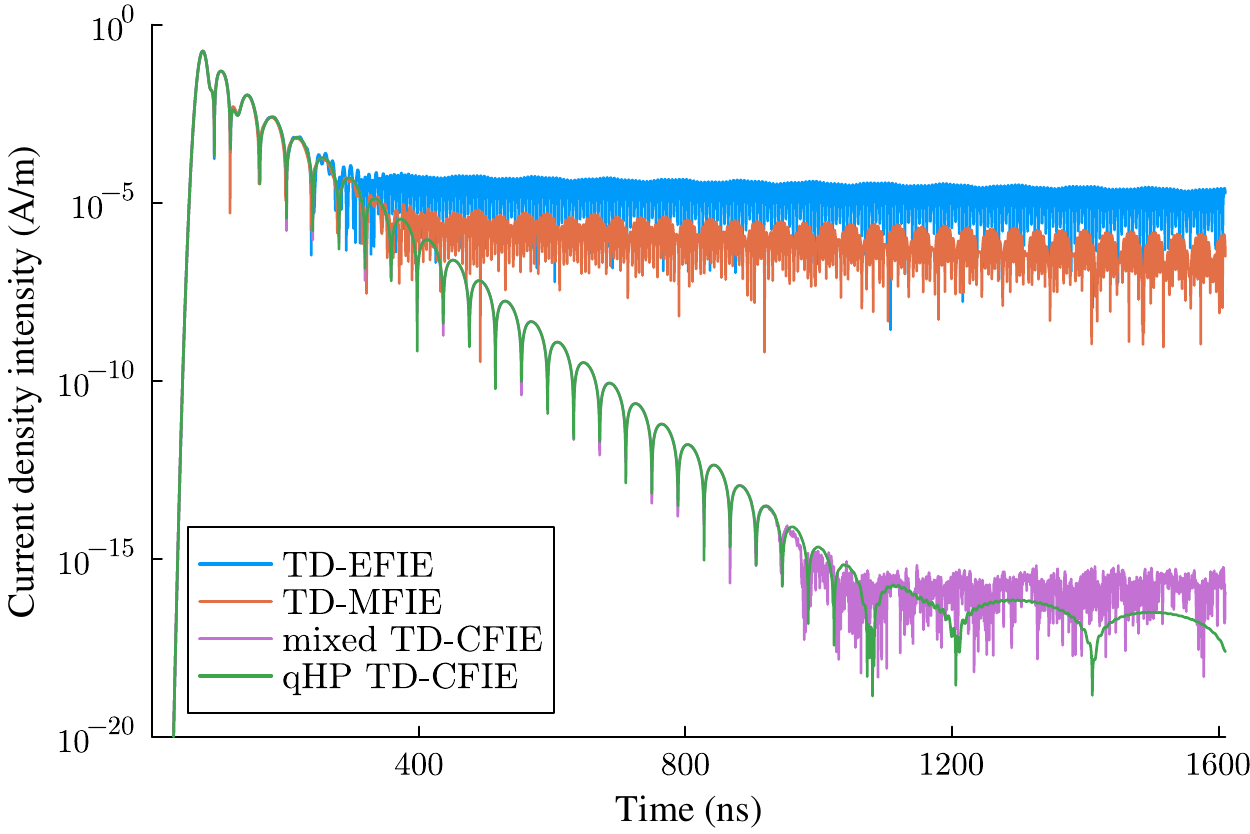}%
\label{fig:current_torus}}
\caption{Intensity of the current density: (a) at the point $\rb = (-0.534, -0.523, -0.644)\mathrm{m}$ on a sphere with radius $1\mathrm{m}$, discretized by 476 triangles; (b) at the point $\rb = (-2.326, -1.885, -0.692)\mathrm{m}$ on a torus with two radii $3\mathrm{m}$ and $1\mathrm{m}$, approximated by 952 triangles. The time step $\Delta t = 0.333 \mathrm{ns}$. The TD-EFIE and the TD-MFIE suffer from resonant instabilities, whereas the TD-CFIEs do not. However, on the torus, the qHP TD-CFIE exhibits a very slowly oscillating mode occurring below the level of machine precision.}
\label{fig:current}
\end{figure*}

Next, we simulate the scattering by the PEC sphere using a moderate time step $\Delta t = 0.333 \mathrm{ns}$. Accordingly, the Gaussian pulse plane wave with width $w = 8\mathrm{m}$ and time of arrival $t_0 = 80\mathrm{ns}$ is used. The surface current density $\jb(\rb, t)$ is computed using four different formulations, and its intensity at the point $\rb = (-0.534, -0.523, -0.644)\mathrm{m}$ are presented in Fig.~\ref{fig:current_sphere}. Before the resonances come into play (at $150 \mathrm{ns}$, approximately), four equations yield similar results.

\begin{figure*}[!t]
\centering
\subfloat[]{\includegraphics[trim={0.3cm 0.8cm 0 0.5cm}, clip, width=0.32\linewidth]{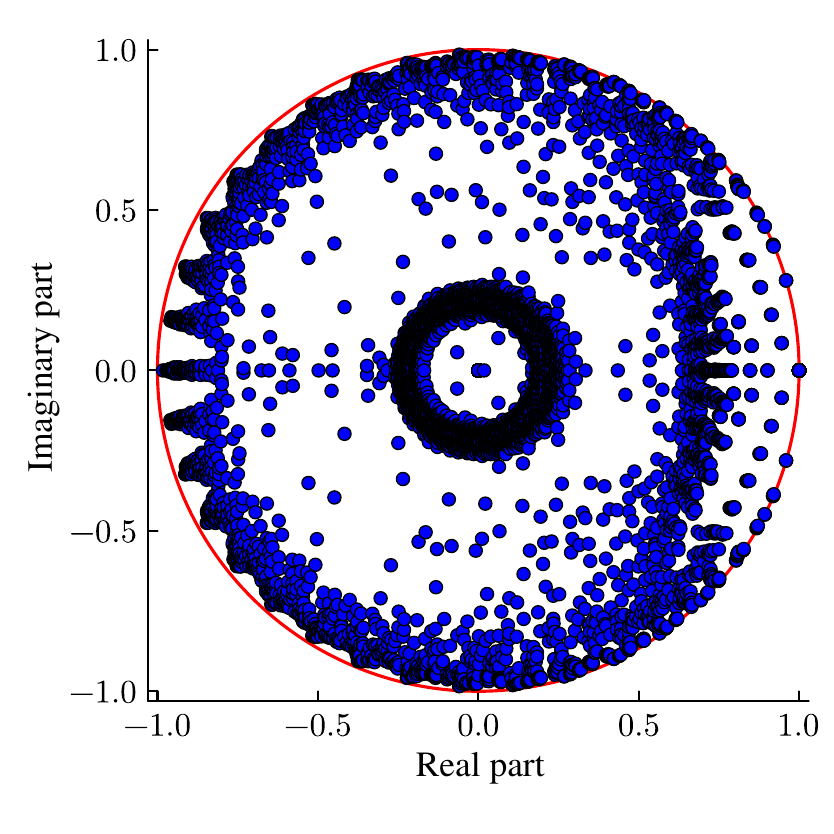}%
\label{fig:polyvals_EFIE_sphere}}
\hfil
\subfloat[]{\includegraphics[trim={0.3cm 0.8cm 0 0.5cm}, clip, width=0.32\linewidth]{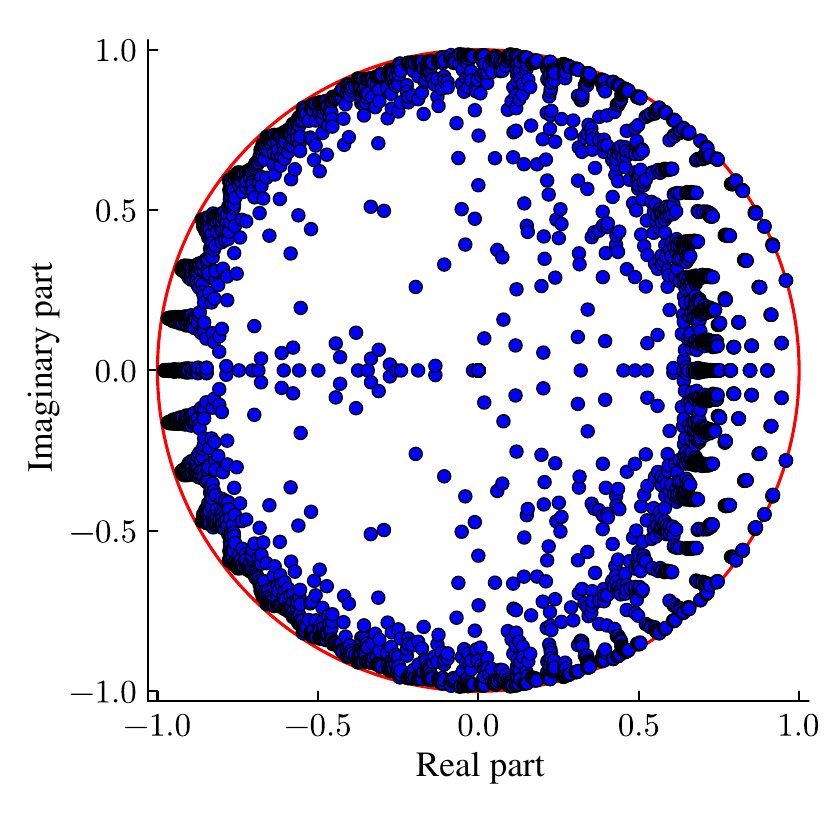}%
\label{fig:polyvals_MFIE_sphere}}
\hfil
\subfloat[]{\includegraphics[trim={0.3cm 0.8cm 0 0.5cm}, clip, width=0.32\linewidth]{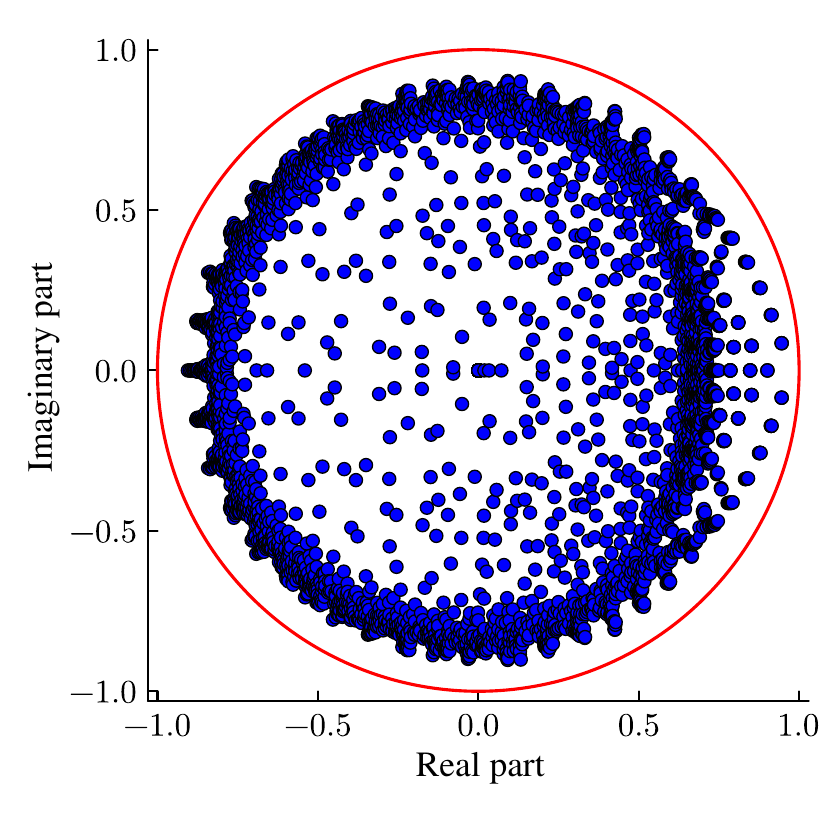}%
\label{fig:polyvals_CFIE_sphere}}
\caption{Polynomial spectrum of (a) the TD-EFIE, (b) the TD-MFIE, and (c) the qHP TD-CFIE for the sphere with radius $1\mathrm{m}$. The time step $\Delta t = 0.333\mathrm{ns}$. The polynomial eigenvalues $1 + 0j$ of the TD-EFIE are corresponding to the constant-in-time solenoidal regime solutions, which are the origin of dc instabilities. The eigenvalues of the form $\lambda = \exp(j\theta)$, with $\theta \neq 0$, residing on the unit circle of the complex plane are corresponding to the interior resonant currents, giving rise to resonant instabilities. The polynomial spectrum of the proposed qHP TD-CFIE is strictly contained in the circle.}
\label{fig:polyvals_sphere}
\end{figure*}

The TD-EFIE and the TD-MFIE are well-known to be prone to spurious resonances. Their polynomial spectrum depicted in Fig.~\ref{fig:polyvals_EFIE_sphere} and Fig.~\ref{fig:polyvals_MFIE_sphere} exhibit the eigenvalues of the form $\lambda = \exp(j\theta)$, with $\theta \neq 0$, which are the root of resonant instabilities. In addition, the eigenvalues around $1 + 0j$ of the TD-EFIE confirm its dc instabilities \cite{ACO+2009}. However, at late times, non-decaying quickly-oscillating errors stemming from the resonances are the most prominent in the solutions.

The mixed TD-CFIE and the qHP TD-CFIE do not suffer from dc and resonant instabilities as the spectrum of the qHP TD-CFIE is strictly contained in the unit circle of the complex plane (see Fig.~\ref{fig:polyvals_CFIE_sphere}). However, the decay of the mixed TD-CFIE solution stops at the level of machine precision ($\approx 10^{-15}$, see Fig.~\ref{fig:current_sphere}). In contrast, the solution of the qHP TD-CFIE \eqref{eq:qHP_TDCFIE} decays infinitely (to $10^{-200}$ and beyond). These results convince us that, by means of separating the solenoidal and irrotational components using the quasi-Helmholtz projectors, then treating them properly, together with explicitly setting vanishing terms to zero, the proposed formulation significantly improves the stability of the solution at late times.

\begin{figure*}[!t]
\centering
\subfloat[]{\includegraphics[width=0.45\linewidth]{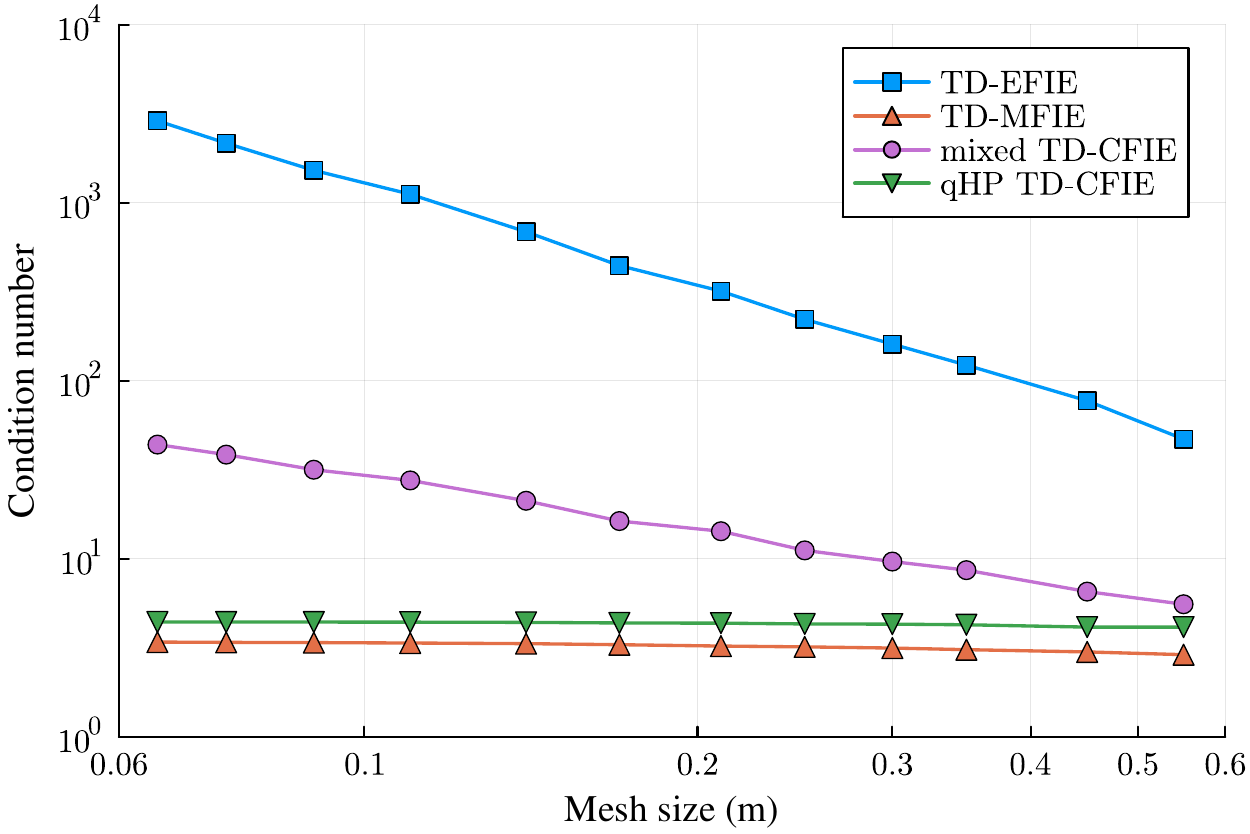}%
\label{fig:cond_h_sphere}}
\hfil \hfil
\subfloat[]{\includegraphics[width=0.45\linewidth]{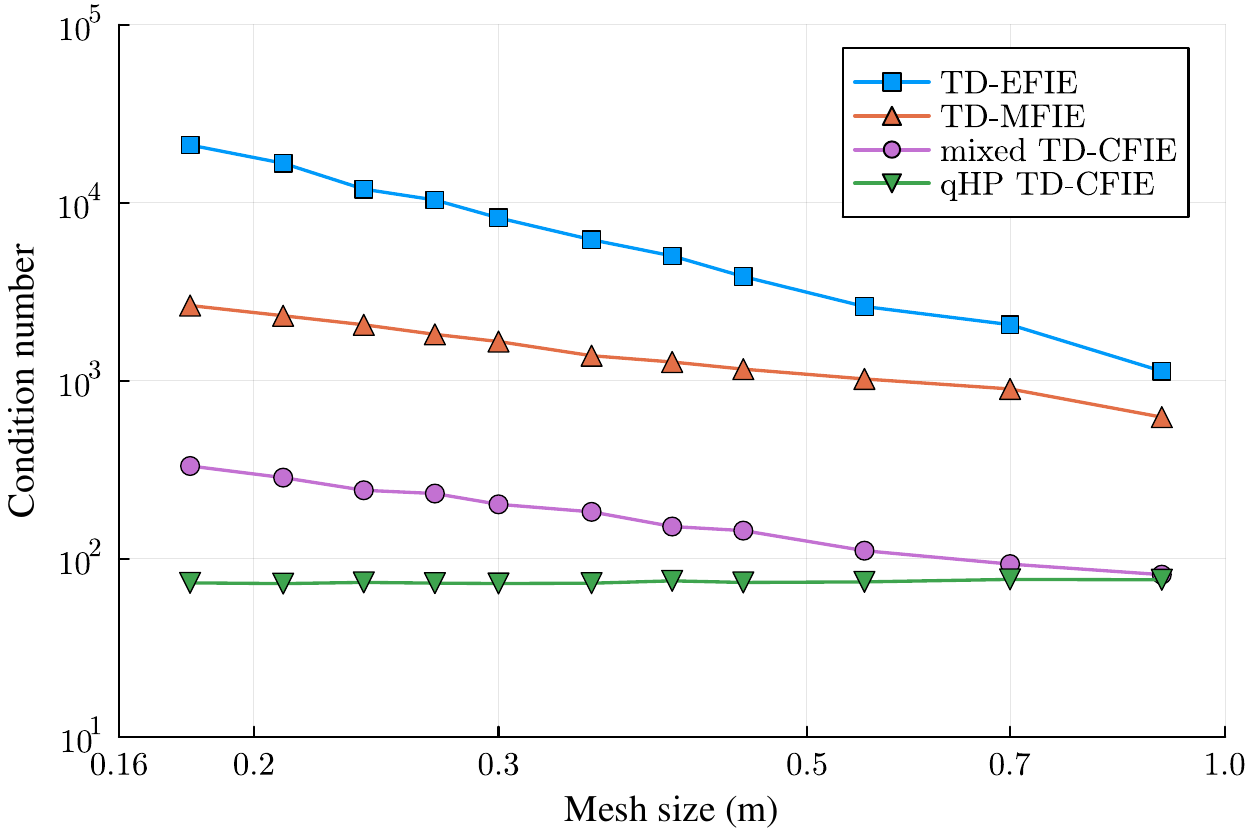}%
\label{fig:cond_h_torus}}
\caption{Condition number of the matrix systems with respect to different mesh sizes: (a) of a sphere with radius $1\mathrm{m}$, while the time step is fixed at $\Delta t = 3.33 \mathrm{ns}$; (b) of a torus with two radii $3\mathrm{m}$ and $1\mathrm{m}$, while the time step is fixed at $\Delta t = 26.667 \mathrm{ns}$. The TD-EFIE and the mixed TD-CFIE become ill-conditioned when the mesh size decreases, whereas the TD-MFIE on the sphere and the proposed qHP TD-CFIE stay well-conditioned. On toroidal surfaces, the TD-MFIE supports the non-trivial kernel of the static outer MFIE operator when the time step is sufficiently large.}
\label{fig:cond_h}
\end{figure*}

\begin{figure*}[!t]
\centering
\subfloat[]{\includegraphics[width=0.45\linewidth]{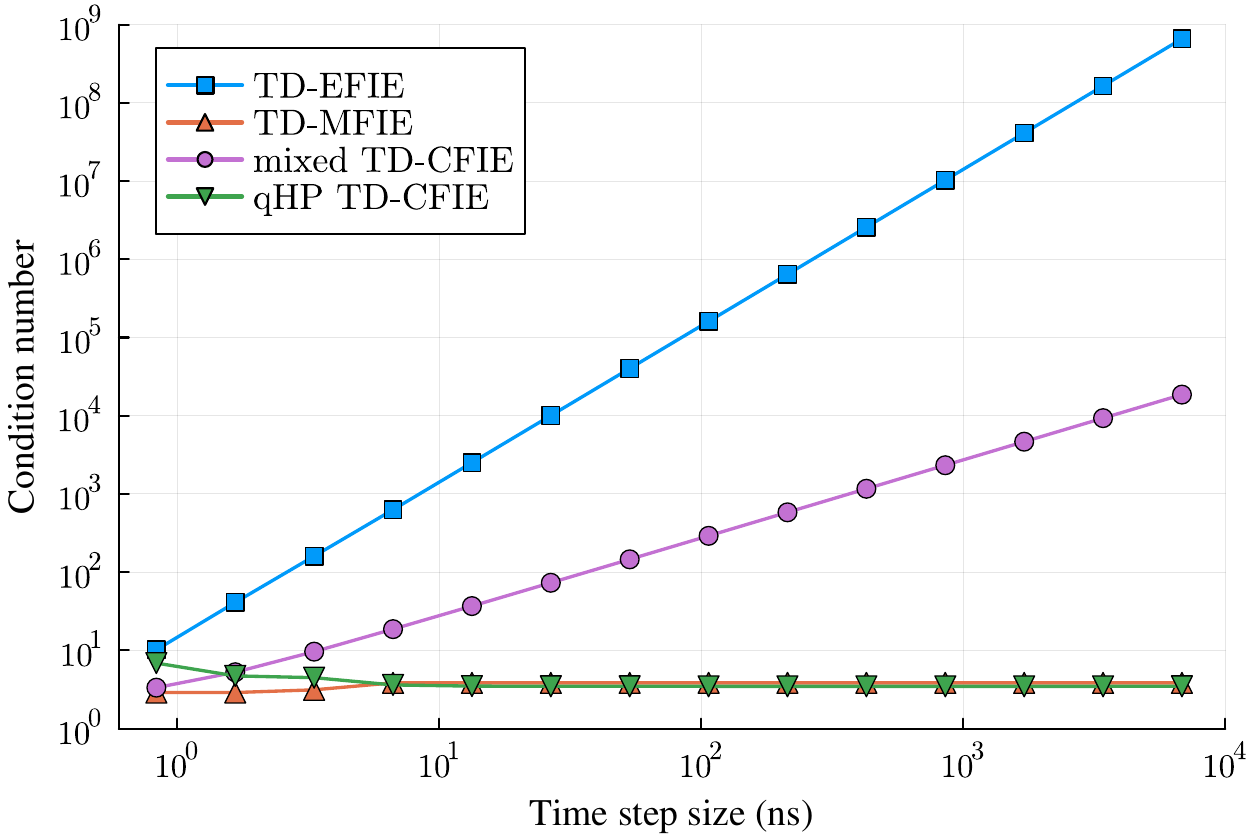}%
\label{fig:cond_dt_sphere}}
\hfil \hfil
\subfloat[]{\includegraphics[width=0.45\linewidth]{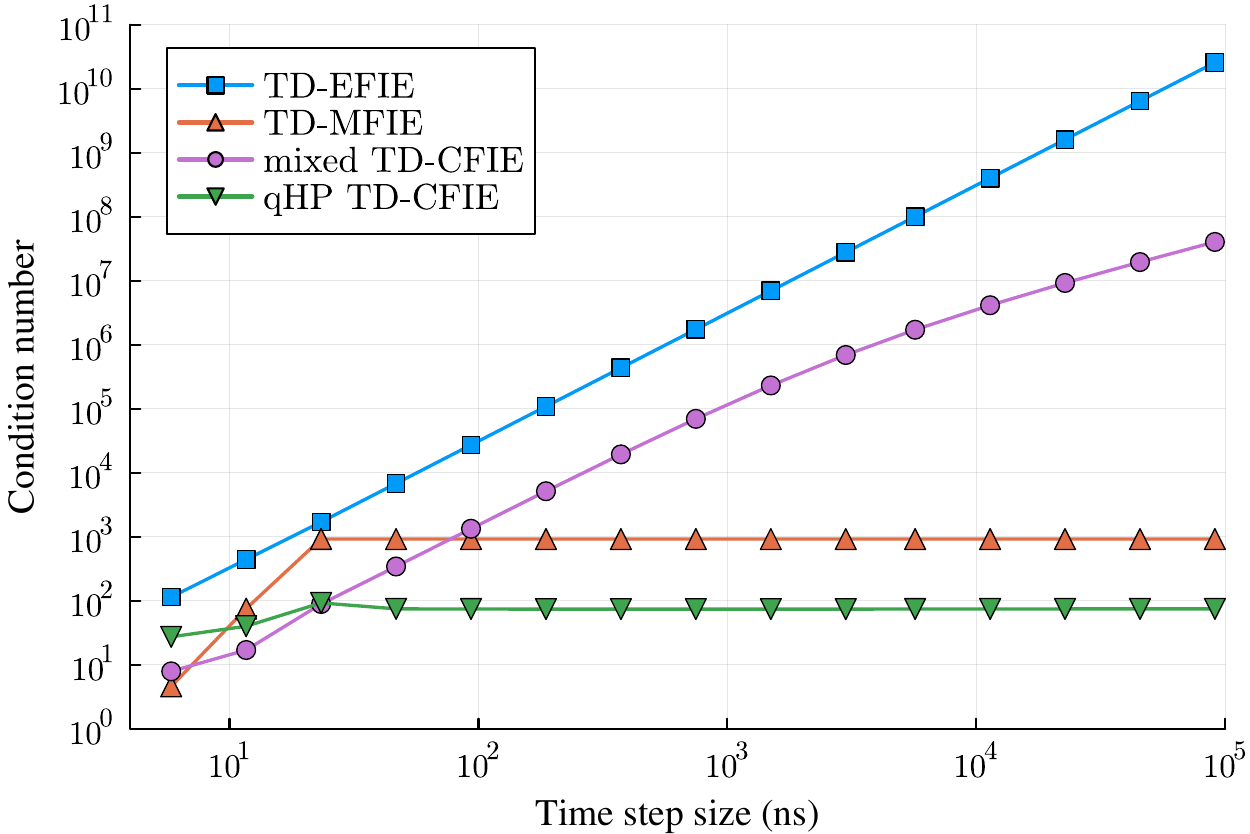}%
\label{fig:cond_dt_torus}}
\caption{Condition number of the matrix systems with respect to different time steps: (a) for a sphere with radius $1\mathrm{m}$, discretized by 476 triangles; (b) for a torus with two radii $3\mathrm{m}$ and $1\mathrm{m}$, discretized by 952 triangles. The TD-EFIE and the mixed TD-CFIE exhibit the large-time step (low-frequency) breakdown, whereas the TD-MFIE and the proposed TD-CFIE do not. However, on toroidal surfaces, the TD-MFIE is affected by the non-trivial nullspaces of the static MFIE operators.}
\label{fig:cond_dt}
\end{figure*}

Besides the computation of numerical solutions, the well-conditioning of the proposed TD-CFIE is also verified. To that end, the condition number of the matrix systems are computed with respect to different time steps and mesh sizes. Fig.~\ref{fig:cond_h_sphere} presents the condition number of the resulting linear systems when the mesh size varies from $0.065\mathrm{m}$ to $0.55\mathrm{m}$ and the time step is fixed at $\Delta t = 3.333 \mathrm{ns}$. In Fig.~\ref{fig:cond_dt_sphere}, the average mesh size is fixed at $0.3 \mathrm{m}$, while the time step varies from $0.833 \mathrm{ns}$ to $6826.7 \mathrm{ns}$. The numerical results show that the condition number of the TD-EFIE and the mixed TD-CFIE increase when either the mesh size decreases, or the time step increases. More importantly, they confirm that the TD-MFIE and the qHP TD-CFIE \eqref{eq:qHP_TDCFIE} do not suffer from both dense discretization breakdown and large-time step breakdown.

\subsection{Torus}

The second experiment concerns the scattering by a PEC torus with large radius $3\mathrm{m}$ and small radius $1\mathrm{m}$ (see Fig.~\ref{fig:geometries}, right). The surface is discretized by 952 triangles with 1428 edges and 476 vertices. The incident plane wave is characterized by the width $w = 8\mathrm{m}$ and the time of arrival $t_0 = 80\mathrm{ns}$, while the considered time step $\Delta t = 0.333 \mathrm{ns}$. The intensity of the surface current density at the point $\rb = (-2.326, -1.885, -0.692)\mathrm{m}$ are shown in Fig.~\ref{fig:current_torus}. Similar to the previous experiment, the TD-EFIE and the TD-MFIE are prone to resonant instabilities. Moreover, on toroidal surfaces, the TD-MFIE is also affected by near-dc instabilities rooted in the nullspaces of the static MFIE operators \cite{CAO+2009a}. This effect can be noticed by the slow oscillation of the amplitude of the TD-MFIE solution. The TD-CFIEs are immune to resonant and dc instabilities. Therefore, their solutions decay, although slowly. From a very late time ($t > 1000 \mathrm{ns}$), the mixed TD-CFIE's solution exhibits numerical noise, while the solution of the qHP TD-CFIE reveals a very slowly oscillating mode that takes place below the level of machine precision. This mode can be interpreted as the weak coupling of the TD-EFIE's dc instabilities and the TD-MFIE's near-dc instabilities. The coupling occurs because the harmonic components (global loops) exist and reside on the active range of both effects.

\begin{figure*}[!t]
\centering
\subfloat[]{\includegraphics[trim={0.3cm 0.8cm 0 0.5cm}, clip, width=0.32\linewidth]{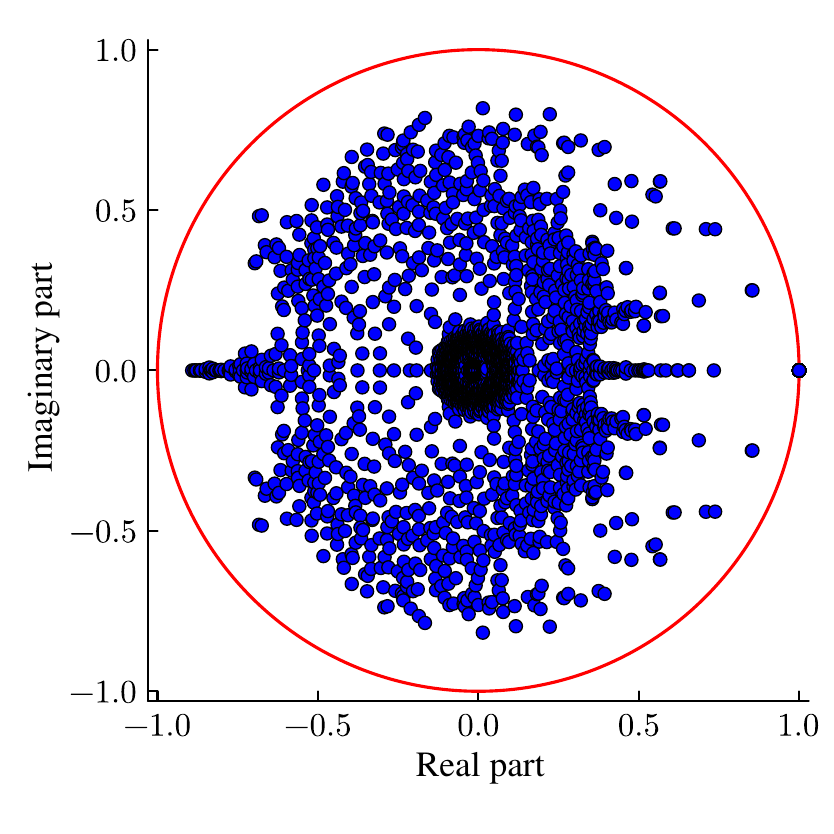}%
\label{fig:polyvals_EFIE_torus}}
\hfil
\subfloat[]{\includegraphics[trim={0.3cm 0.8cm 0 0.5cm}, clip, width=0.32\linewidth]{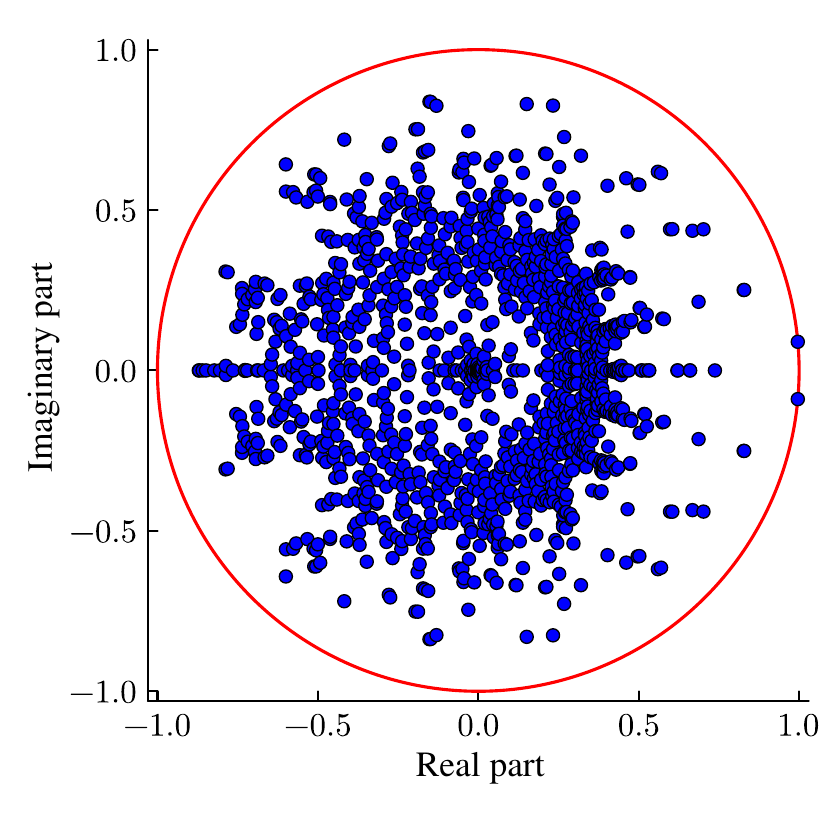}%
\label{fig:polyvals_MFIE_torus}}
\hfil
\subfloat[]{\includegraphics[trim={0.3cm 0.8cm 0 0.5cm}, clip, width=0.32\linewidth]{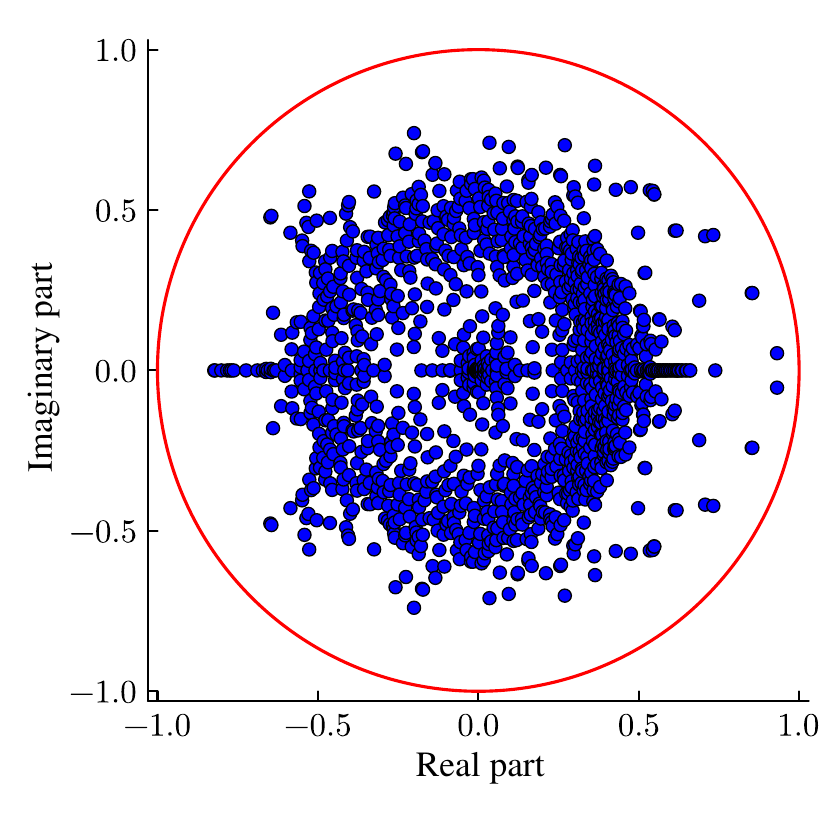}%
\label{fig:polyvals_CFIE_torus}}
\caption{Polynomial spectrum of (a) the TD-EFIE, (b) the TD-MFIE, and (c) the qHP TD-CFIE for the torus with two radii $3\mathrm{m}$ and $1\mathrm{m}$. The time step $\Delta t = 3.333\mathrm{ns}$. The TD-EFIE exhibits dc instabilities associated with the eigenvalues $1 + 0j$. The TD-MFIE has two conjugate eigenvalues residing on the unit circle, which are associated with instabilities rooted in the nullspaces of the static MFIE operators. The polynomial spectrum of the qHP TD-CFIE is strictly contained in the circle. However, it also has two poles approaching $1 + 0j$ from inside when the time step and the mesh size decrease, revealing exponentially-decaying near-dc instabilities of the proposed TD-CFIE.}
\label{fig:polyvals_torus}
\end{figure*}

In order to gain a deeper insight into the spectral properties of the TD-EFIE, the TD-MFIE and the qHP TD-CFIE, their polynomial eigenvalue analysis are performed and illustrated in Fig.~\ref{fig:polyvals_torus}. The time step $\Delta t = 3.333\mathrm{ns}$ (ten times larger than the one considered in Fig.~\ref{fig:current_torus}). The eigenvalues associated with the resonant currents are not present since the frequency range we can model is below the lowest resonant frequency. However, the eigenvalues $1 + 0j$ are present in the spectrum of the TD-EFIE, which are the origin of dc instabilities. In accordance with \cite{CAO+2009a}, the TD-MFIE has two conjugate polynomial eigenvalues residing on the unit circle of the complex plane, which are associated with non-decaying near-dc instabilities. In the spectrum of the qHP TD-CFIE, these two eigenvalues are shifted inside the circle (see Fig.~\ref{fig:polyvals_CFIE_torus}), implying the stability of the proposed formulation in typical applications. Nevertheless, when the temporal and spatial mesh density increases, these two poles tend toward $1 + 0j$ from inside, revealing near-dc instabilities corresponding to exponentially-decaying error components. These instabilities are rooted in the weak coupling of the TD-EFIE's dc instabilities and the MFIE's near-dc instabilities.

Figures~\ref{fig:cond_h_torus} and \ref{fig:cond_dt_torus} again corroborate the well-conditioning of the proposed TD-CFIE \eqref{eq:qHP_TDCFIE} when the spatial mesh is fine or the time step is large. Please note that if $c\Delta t \ge D$ ($D$ is the diameter of the scatterer), the conditioning behavior of the TD-MFIE coincides with that of the static outer MFIE operator. This fact explains the increase of the condition number of the TD-MFIE along the resolution of the spatial mesh in Fig.~\ref{fig:cond_h_torus}. It can also explain the behavior in Fig.~\ref{fig:cond_dt_torus}, where the condition number of the TD-MFIE firstly increases, and afterwards it stays constant when the time step is sufficiently large.

\section{Conclusion}
\label{sec:conclusion}

The novel TD-CFIE formulation is proposed by means of the quasi-Helmholtz projectors. The solenoidal and irrotational components are separated, then they are properly rescaled by appropriate factors. The stabilized Calder\'{o}n preconditioned TD-EFIE and the stabilized symmetrized TD-MFIE are linearly combined in a Calder\'{o}n-like fashion. The semi-discrete TD-CFIE formulation is immune to dc instabilities, resonant instabilities, dense discretization breakdown and large-time step breakdown. The proposed qHP TD-CFIE is then temporally discretized using a compatible pair of temporal basis functions. The MOT scheme yields a stable and accurate solution. The formulation developed in this paper can be applied to both simply-connected and multiply-connected geometries.

Unfortunately, near-dc instabilities stemming from the weak coupling of the TD-EFIE's dc instabilities and the TD-MFIE's near-dc instabilities may prevent us from obtaining extremely accurate solutions for multiply-connected surfaces. In the forthcoming study, another rescaling procedure should be introduced to eliminate the non-trivial kernel of the TD-EFIE (for instance based on the procedures proposed in \cite{BCA2015} and \cite{BCA2015b}). If successful, the proposed formulation will resolve all numerical instabilities and breakdowns of time-domain integral equations at moderate to extremely low frequencies.

% or


% \appendices

% % you can choose not to have a title for an appendix
% % if you want by leaving the argument blank
% \section{}

% \section*{Acknowledgment}

\ifCLASSOPTIONcaptionsoff
  \newpage
\fi

% \IEEEtriggeratref{36}
% \IEEEtriggercmd{\enlargethispage{-5in}}


\bibliographystyle{IEEEtran}
\bibliography{abrv_ref.bib}


\begin{IEEEbiography}[{\includegraphics[width=1in,height=1.25in,clip,keepaspectratio]{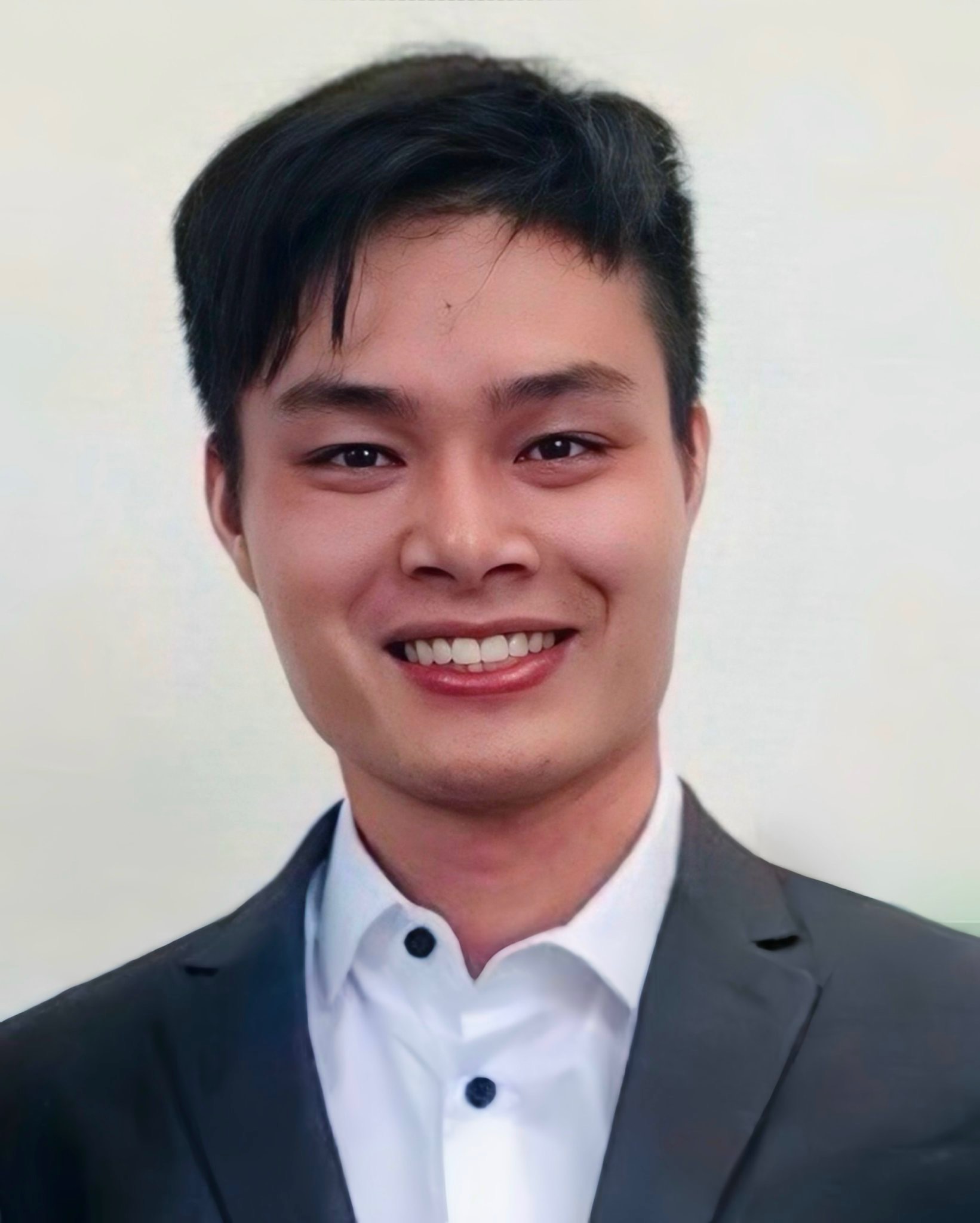}}]{Van Chien Le}
received the M.Sc. degree in applied mathematics from Hanoi University of Science and Technology, Hanoi, Vietnam in 2018, and the Ph.D. degree in mathematics from Ghent University, Ghent, Belgium in 2022.

Since 2022, he has been a postdoctoral researcher at the Department of Information and Technology, Ghent University, Ghent, Belgium. His current research interests include numerical analysis of finite element and boundary element methods for electromagnetics, and stable and accurate time-domain boundary integral equation solvers for electromagnetic scattering.

Dr. Le was a recipient of the IEEE Ulrich L. Rohde innovative conference paper award on computational techniques in electromagnetics at the International Conference on Electromagnetics in Advanced Applications in 2023. He also co-authored a paper receiving an honorable mention at the IEEE International Symposium on Antennas and Propagation and USNC-URSI Radio Science Meeting in 2023.
\end{IEEEbiography}

\begin{IEEEbiography}[{\includegraphics[width=1in,height=1.25in,clip,keepaspectratio]{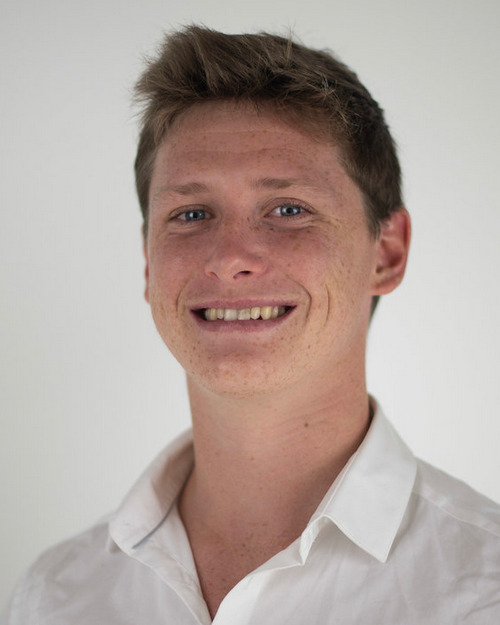}}]{Pierrick Cordel}
(Graduate Student Member, IEEE) received the M.Sc.Eng. degree from the Ecole Nationale Sup\'{e}rieure des Mines de Nancy, France, in 2021, as well as the M.Sc. degree in industrial mathematics from the University of Luxembourg, Luxembourg, in 2021. 

He is currently doing the Ph.D. thesis with the Institute Politecnico di Torino, Italy. His research focuses on preconditioned and fast solution of boundary element methods frequency domain and time-domain integral equations.
\end{IEEEbiography}

\begin{IEEEbiography}[{\includegraphics[width=1in,height=1.25in,clip,keepaspectratio]{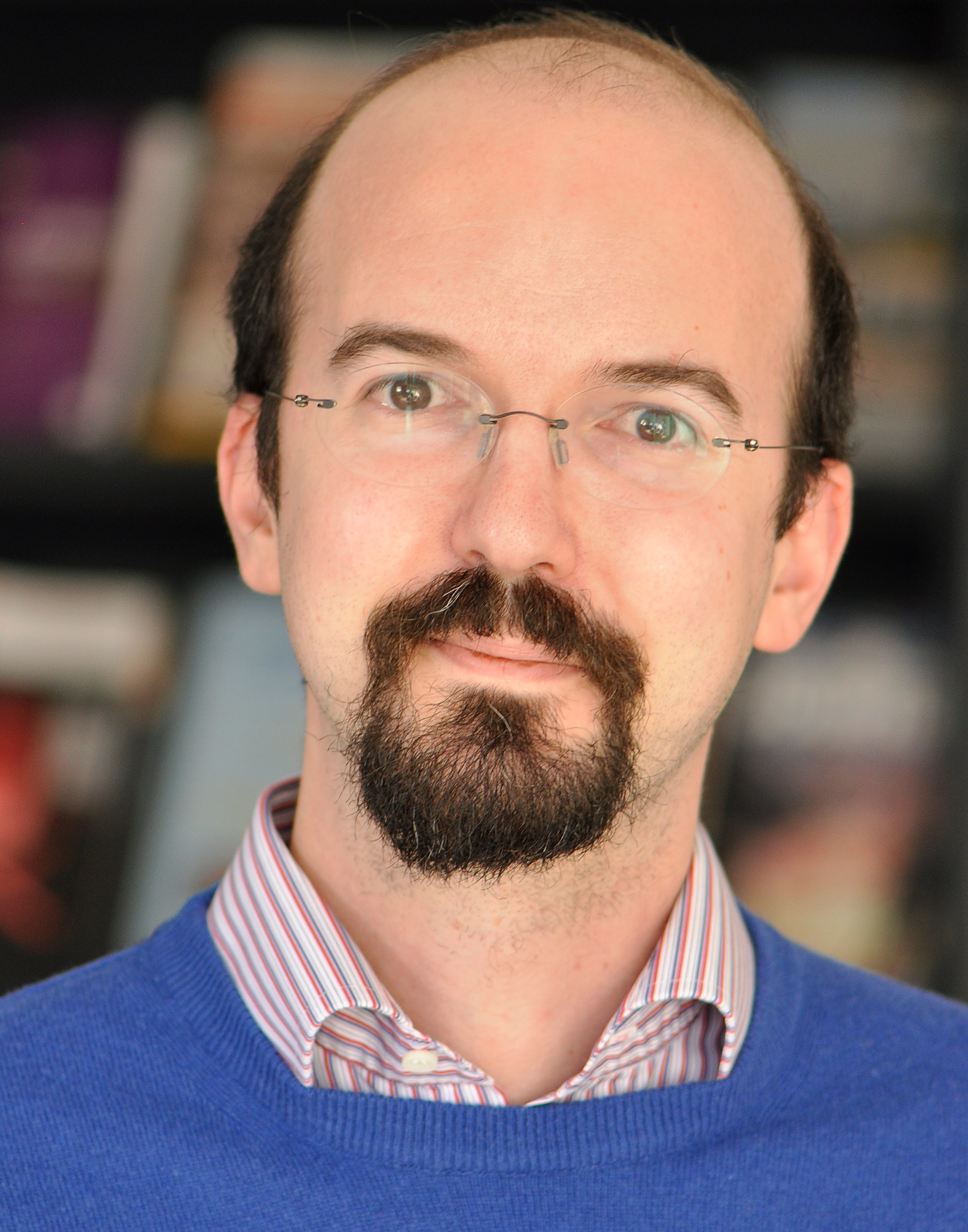}}]{Francesco P. Andriulli}
(S'05--M'09--SM'11--F'23) received the Laurea in electrical engineering from the Politecnico di Torino, Italy, in 2004, the MSc in electrical engineering and computer science from the University of Illinois at Chicago in 2004, and the PhD in electrical engineering from the University of Michigan at Ann Arbor in 2008. From 2008 to 2010 he was a Research Associate with the Politecnico di Torino. From 2010 to 2017 he was an Associate Professor (2010-2014) and then Full Professor with the École Nationale Supérieure Mines-Télécom Atlantique (IMT Atlantique, previously ENST Bretagne), Brest, France. Since 2017 he has been a Full Professor with the Politecnico di Torino, Turin, Italy. His research interests are in computational electromagnetics with focus on frequency- and time-domain integral equation solvers, well-conditioned formulations, fast solvers, low-frequency electromagnetic analyses, and modeling techniques for antennas, wireless components, microwave circuits, and biomedical applications with a special focus on brain imaging.
 
Prof. Andriulli received several best paper awards at conferences and symposia (URSI NA 2007, IEEE AP-S 2008, ICEAA IEEE-APWC 2015) also in co-authorship with his students and collaborators (ICEAA IEEE-APWC 2021, EMTS 2016, URSI-DE Meeting 2014, ICEAA 2009) with whom received also a second prize conference paper (URSI GASS 2014), a third prize conference paper (IEEE–APS 2018), seven honorable mention conference papers (ICEAA 2011, URSI/IEEE–APS 2013, 4 in URSI/IEEE–APS 2022, URSI/IEEE–APS 2023) and other three finalist conference papers (URSI/IEEE-APS 2012, URSI/IEEE-APS 2007, URSI/IEEE-APS 2006, URSI/IEEE–APS 2022). A Fellow of the IEEE, he is also the recipient of the 2014 IEEE AP-S Donald G. Dudley Jr. Undergraduate Teaching Award, of the triennium 2014-2016 URSI Issac Koga Gold Medal, and of the 2015 L. B. Felsen Award for Excellence in Electrodynamics.
 
Prof. Andriulli is a member of Eta Kappa Nu, Tau Beta Pi, Phi Kappa Phi, and of the International Union of Radio Science (URSI). He is the Editor-in-Chief of the \textit{IEEE Antennas and Propagation Magazine} and he serves as a Track Editor for the \textit{IEEE Transactions on Antennas and Propagation}, and as an Associate Editor for \textit{URSI Radio Science Letters}. He served as an Associate Editor for the \textit{IEEE Transactions on Antennas and Propagation}, \textit{IEEE Antennas and Wireless Propagation Letters}, \textit{IEEE Access and IET-MAP}. He is the PI of the ERC Consolidator Grant: \textit{321 – From Cubic$^3$ To$^2$ Linear$^1$ Complexity in Computational Electromagnetics.}
\end{IEEEbiography}
\newpage

\begin{IEEEbiography}
[{\includegraphics[width=1in,height=1.25in,clip,keepaspectratio]{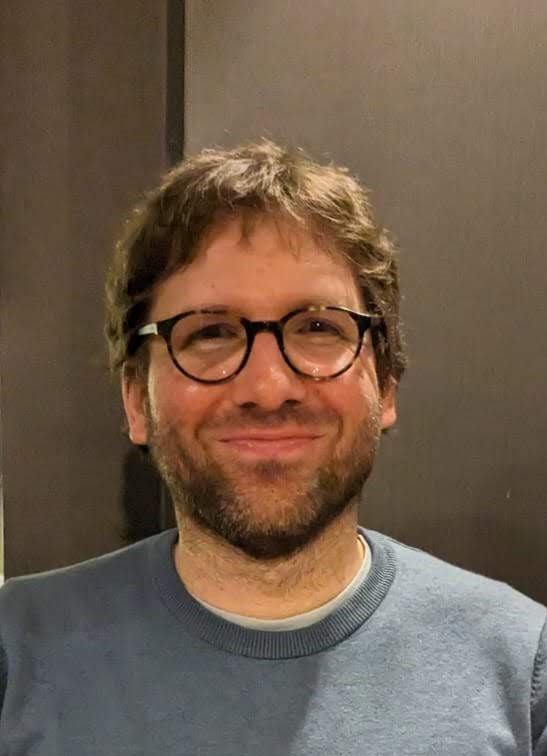}}]{Kristof Cools} is full professor with the Department of Information Technology at Ghent University. Before joining Ghent University he has held positions at the TU Delft and the University of Nottingham. Kristof Cools has been awarded the ICEAA Young Scientist Award 2009, the URSI Young Scientist Award 2014, and the Rohde Most Innovative Conference Paper Award 2023. He has authored over 33 papers published in peer reviewed international journals and over 120 conference contributions, resulting in a Google Scholar recorded h-index of 19 and over 1800 citations. Dr Cools has received funding from institutional, regional, industrial, national, and European funding bodies, including a Marie Skłodowska-Curie career integration grant and ERC consolidator project Boundary Elements for Time-Domain Domain Decomposition Methods. 
\end{IEEEbiography}


%\newpage



%\vfill


\end{document}